\documentclass[10pt,titlepage]{amsart} 
\usepackage{bookman,amssymb,mathrsfs}
\usepackage[author-year]{amsrefs}%
   \date{September 20, 2004}
   \title[Dynamical Walks]{Exceptional Times and 
   Invariance\\ for dynamical Random Walks 
   } 
   \thanks{The research of 
      D.\@ Kh.\@ is partially supported by a grant from the NSF} 
   \address{Department\@ of Mathematics\\The University\@ of Utah\\ 
      155 S.\@ 1400 E.\\Salt Lake City, UT 84112--0090} 
   \author[D. Khoshnevisan]{Davar Khoshnevisan} 
   \email{davar@math.utah.edu} 
   \urladdr{http://www.math.utah.edu/\~{}davar} 
   \author[D.\@ A.\@ Levin]{David A.\@ Levin} 
   \address{Department\@ of Mathematics\\The University\@ of Utah\\ 
      155 S.\@ 1400 E.\\Salt Lake City, UT 84112--0090} 
   \email{levin@math.utah.edu} 
   \urladdr{http://www.math.utah.edu/\~{}levin} 
   \author[P.\@ M\'endez]{Pedro J. M\'endez-Hern\'andez} 
   \address{Department\@ of Mathematics\\The University\@ of Utah\\ 
      155 S.\@ 1400 E.\\Salt Lake City, UT 84112--0090}
   \curraddr{Escuela de Matematica\\Universidad de Costa Rica\\
      San Pedro de Montes de Oca\\Costa Rica}
   \email{mendez@math.utah.edu,pmendez@emate.ucr.ac.cr}
   \urladdr{http://www.math.utah.edu/\~{}mendez} 
 
\theoremstyle{plain}{ 
\newtheorem{theorem}{Theorem}[section]} 
\theoremstyle{plain}{ 
\newtheorem{proposition}[theorem]{Proposition}} 
\theoremstyle{plain}{ 
   \newtheorem{lemma}[theorem]{Lemma}} 
\theoremstyle{plain}{ 
   \newtheorem{pbm}[theorem]{Problem}} 
\theoremstyle{plain}{ 
   \newtheorem{corollary}[theorem]{Corollary}} 
\theoremstyle{definition}{ 
   \newtheorem{definition}[theorem]{Definition}} 
\theoremstyle{definition}{ 
   \newtheorem{example}[theorem]{Example}} 
\theoremstyle{remark}{ 
   \newtheorem{remark}[theorem]{Remark}} 
\numberwithin{equation}{section} 
\newcommand{\K}{\mathrm{K}} 
\newcommand{\M}{\mathrm{M}} 
\newcommand{\Hdim}{\dim_{_{\mathscr H}}} 
\newcommand{\Kdim}{\overline{\dim}_{_{\mathscr M}}} 
 
\newcommand{\Pdim}{\dim_{_{\mathscr P}}} 
\newcommand{\F}{\mathscr{F}}
 
\newcommand{\e}{\varepsilon} 
\newcommand{\s}{\sigma} 
\renewcommand{\d}{\delta} 
\newcommand{\PN}{\P_{\mathscr{N}}} 
\newcommand{\EN}{\E_{\mathscr{N}}} 
\renewcommand{\P}{\mathrm{P}} 
\newcommand{\ee}{\mathfrak{e}} 
\newcommand{\E}{\mathrm{E}} 
\newcommand{\R}{\mathbf{R}} 
\newcommand{\fclass}{\pmb{\mathscr{C}}_{bz}} 
\newcommand{\inprob}{\stackrel{{\P}}{\to}} 
\newcommand{\floor}[1]{\lfloor #1 \rfloor} 

\begin{document} 
\begin{abstract} 
   Consider a sequence $\{X_i(0)\}_{i=1}^n$ 
   of i.i.d.\@ random variables. Associate 
   to each $X_i(0)$ an independent mean-one 
   Poisson clock. Every time a clock rings 
   replace that $X$-variable by an independent 
   copy and restart the clock. In this way, 
   we obtain i.i.d.\@ stationary processes 
   $\{X_i(t)\}_{t\ge 0}$ 
   ($i=1,2,\cdots$) whose invariant distribution is 
   the law $\nu$ of $X_1(0)$. 
 
   \ocite{benjamini} introduced the dynamical 
   walk $S_n(t)=X_1(t)+\cdots+X_n(t)$, and proved 
   among other things that the LIL holds for $n\mapsto S_n(t)$ 
   for all $t$. In other words, the LIL is dynamically 
   stable. Subsequently~\ycite{KLM}, we showed that 
   in the case that the $X_i(0)$'s are standard normal, 
   the classical integral test is not dynamically stable. 
 
   Presently, we study the set of times $t$ 
   when $n\mapsto S_n(t)$ 
   exceeds a given envelope infinitely often. Our analysis is 
   made possible thanks to a connection to the 
   Kolmogorov $\e$-entropy. When used in conjunction 
   with the invariance principle of this paper, this connection 
   has other interesting by-products some of which we relate. 
 
   We prove also that the 
   infinite-dimensional process 
   $t\mapsto S_{\lfloor n\bullet\rfloor}(t)/\sqrt n$ 
   converges weakly in $\mathscr{D}
   (\mathscr{D}([0,1]))$ to the 
   Ornstein--Uhlenbeck process 
   in $\mathscr{C}([0,1])$. For this we assume only 
   that the increments have mean zero and variance one.
 
   In addition, we extend a result of~\ocite{benjamini} 
   by proving that if the $X_i(0)$'s are lattice, mean-zero 
   variance-one, and possess $2+\e$ finite absolute 
   moments for some $\e>0$, 
   then the recurrence of the origin is dynamically stable. 
   To prove this we derive a gambler's 
   ruin estimate that is valid for all lattice random walks 
   that have mean zero and finite variance. We believe 
   the latter may be of independent interest. 
\end{abstract} 
\maketitle\tableofcontents 
\section{Introduction and Main Results} 
 
Let $\{\xi_j^k\}_{j,k=0}^\infty$ denote a double-array 
of i.i.d.\@ real-valued random variables with common 
distribution $\nu$. Also let 
$\{\mathscr{C}_n\}_{n=1}^\infty$ denote a sequence of rate-one Poisson clocks 
that are totally independent from themselves as well as the 
$\xi$'s. If the jump times of $\mathscr{C}_n$ are denoted by 
$0=\tau_n(0)<\tau_n(1)<\tau_n(2)<\cdots$, then we define the discrete-time 
function-valued process $X =\{ X_n(t);\, t\ge 0\}_{n=1}^\infty$ 
as follows: For all $n\ge 1$, 
\begin{equation} 
   X_n (t)  =   \xi_n^k,\quad 
   \text{if }\tau_n(k)\le t < \tau_n(k+1). 
\end{equation} 
For every $n\ge 1$, $X_n$ is the random step function which 
starts, at time zero, at the value $\xi^0_n$. 
Then it proceeds iteratively 
by replacing its previous value by an independent copy 
every time the clock $\mathscr{C}_n$ rings. As a process indexed by 
$t$, $t\mapsto (X_1(t),X_2(t),\ldots)$ is a stationary Markov process 
in $\R^\infty$, and its invariant measure is $\nu^\infty$. 
 
The \emph{dynamical walk} corresponding to the $X$'s is the random field 
\begin{equation} 
   S_n(t) = X_1(t)+\cdots+X_n(t),\qquad 
   t\ge 0,\ n\ge 1. 
\end{equation} 
One can think of the ensuing object in different ways. 
We take the following points of view interchangeably: 
\begin{enumerate} 
   \item For a given $t\ge 0$, $\{S_n(t)\}_{n=1}^\infty$ 
      is a classical random 
      walk with increment-distribution $\nu$. 
   \item For a given $n\ge 1$, $\{S_n(t)\}_{t\ge 0}$ is a right-continuous 
      stationary Markov process in $\R$ 
      whose invariant measure is $\nu*\cdots*\nu$ ($n$ times). 
   \item The process $\{S_n(\cdot)\}_{n=1}^\infty$ 
      is a random walk with values in the Skorohod space 
      $\mathscr{D}([0,1])$. 
   \item  The $\R^\infty$-valued process 
      $t\mapsto (S_1(t),S_2(t),\ldots)$ is right-continuous, 
      stationary, and Markov. Moreover, 
      its invariant measure is the evolution law of a 
      classical random walk with increment-distribution $\nu$. 
\end{enumerate} 
 
Dynamical walks were introduced 
recently by I.~Benjamini, O.~H\"aggstr\"om, 
Y.~Peres, and J.~Steif~\ycite{benjamini} who posed the following question: 
\begin{equation}\label{DS}\begin{split} 
   &\text{Which a.s.-properties of the classical random walk 
      [$\nu$]}\\[-1mm] 
   &\text{hold simultaneously for all $t\in [0,1]$?} 
\end{split}\end{equation} 
 
Random-walk properties that satisfy~\eqref{DS} are called 
\emph{dynamically stable}; all others 
are called \emph{dynamically sensitive}. This definition was introduced by 
\ocite{benjamini} who proved, among many other things, that 
if $\int_{-\infty}^\infty x^2\, \nu(dx)$ is finite 
then: 
\begin{equation}\label{eq:LIL-DS} 
   \text{The law of the iterated logarithm is dynamically stable.} 
\end{equation} 
In order to write this out properly, let us assume, 
without loss of generality, that 
$\int x\,\nu(dx)=0$ and $\int x^2\,\nu(dx)=1$. 
For any non-decreasing measurable function $H:\R_+\to\R_+$ define 
\begin{equation}\label{eq:KLM} 
   \Lambda_H =  \left\{ t\in[0,1]:\ 
   S_n(t) > H(n)\sqrt{n} 
   \text{ infinitely often }[n] \right\}. 
\end{equation} 
In words, $\Lambda_H$ denotes the set of times 
$t\in[0,1]$ when $H(n)\sqrt{n}$ fails to be 
in the upper class [in the sense of 
P.\@ L\'evy] of the process $\{ S_n (t)\}_{n=1}^\infty$. 
According to the Hewitt--Savage zero-one law, the event 
$\{\Lambda_H \neq \varnothing\}$ has probability zero or one. 
 
Now set $H(n) = \sqrt{2c\ln\ln n}$. 
Then, dynamical stability of the LIL (\ref{eq:LIL-DS}) is 
equivalent to the statement that 
$\Lambda_H=\varnothing$ a.s.~if $c>1$, whereas 
$\Lambda_H=[0,1]$ a.s.~if $c<1$. 
After ignoring a null set, we can write this in the following 
more conventional form: 
\begin{equation} 
   \limsup_{n\to\infty} \frac{S_n(t)}{\sqrt{2n\ln\ln n}} =1, 
   \qquad{}^\forall t\in[0,1]. 
\end{equation} 
Despite this, in the case that $\nu$ is standard 
normal, we have: 
\begin{equation}\label{eq:UC-Ds}\begin{split} 
   &\text{The characterization of the upper class of a Gaussian}\\[-1mm] 
   &\text{random walk is dynamically sensitive.} 
\end{split}\end{equation} 
Let $\Phi$ denote the standard 
normal distribution function and define 
$\bar\Phi=1-\Phi$. Recall 
that $H(n)\sqrt{n}$ is in the upper class of $S_n(0)$ 
if and only if 
$\int_1^\infty H^2(t)\bar\Phi(H(t))\, dt/t<\infty$~\cite{erdos}. 
Then, \eqref{eq:UC-Ds} is a consequence of Erd\H{o}s's theorem, 
used in conjunction with 
the following recent result~\cite{KLM}*{Theorem 1.5}: 
\begin{equation}\label{KLM:T1.5} 
   \Lambda_H \neq \varnothing \ \Longleftrightarrow\ 
   \int_1^\infty H^4(t) 
   \bar\Phi(H(t)) \, \frac{dt}{t}=\infty. 
\end{equation} 
 
This leaves open the following natural question: Given a non-decreasing 
function $H$, how large is the collection of all times $t\in[0,1]$ at which 
$H$ fails to be in the upper class of $\{S_n(t)\}_{n=1}^\infty$? 
In other words, we ask, ``\emph{How large is $\Lambda_H$}''? 
Define 
\begin{equation}\label{eq:delta} 
   \d (H) =  \sup\left\{ \zeta>0:\ 
   \int_1^\infty H^\zeta(t)\bar\Phi(H(t))\, 
   \frac{dt}{t} <\infty\right\}, 
\end{equation} 
where $\sup\varnothing= 0$. 
The following 
describes the size of $\Lambda_H$ in terms 
of its Hausdorff--Besicovitch dimension $\Hdim\Lambda_H$. 
 
\begin{theorem}\label{thm:Hdim} 
   Suppose $\nu$ is standard normal and $H:\R_+\to\R_+$ is non-random and 
   non-decreasing. Then with probability one, 
   \begin{equation}\label{eq:Hdim} 
      \Hdim \Lambda_H 
      =\min \left( 1 ~,~ \frac{4-\d (H)}{2}\right), 
   \end{equation} 
   where $\Hdim A<0$ means that $A$ is empty. 
\end{theorem} 
In order to prove this we develop a series of technical 
results of independent interest. We describe one of them next. 
 
First define, for any $\e>0$, 
$k= \K_E(\e)$ to be the maximal number of points 
$x_1,\ldots,x_k\in E$ such that whenever $i\neq j$, 
$|x_i-x_j|\ge\e$. The function $\K_E$ 
is known as the \emph{Kolmogorov $\e$-entropy} of 
$E$~\cite{tihomirov}, as well 
as the \emph{packing number (or function)} of $E$~\cite{Mattila}. 
Now suppose $\{z_j\}_{j=1}^\infty$ is any sequence of real 
numbers that satisfies 
\begin{equation}\label{eq:z} 
   \inf_n z_n \ge 1,\ 
   \lim_{n\to\infty}z_n=\infty,\ \text{ and }\ 
   \lim_{n\to\infty}\frac{z_n}{n^{1/4}}=0. 
\end{equation} 
Then we have the following estimate; 
it expresses how the geometry of $E$ affects 
probabilities of moderate deviations. 
\begin{theorem}\label{thm:genest} 
   Let $\nu$ be standard normal. 
   In addition, choose and fix a 
   sequence $\{ z_j\}_{j=1}^\infty$ that 
   satisfies \eqref{eq:z}. 
   Then there exists a finite 
   constant $A_{\ref{eq:genest}}>1$ such that for all $n\ge 1$ 
   and all non-empty non-random measurable sets $E\subseteq[0,1]$, 
   \begin{equation}\label{eq:genest} 
      A_{\ref{eq:genest}}^{-1} \K_E\left( 
      \frac{1}{z_n^2}\right) \bar\Phi(z_n) \le 
      \P\left\{ \sup_{t\in E} S_n(t) \ge z_n\sqrt{n} 
      \right\} \le A_{\ref{eq:genest}} 
      \K_E\left( \frac{1}{z_n^2}\right) \bar\Phi(z_n). 
   \end{equation} 
\end{theorem} 
 
Theorem~\ref{thm:inttest2} below 
appeals to Theorem~\ref{thm:genest} to characterize 
all non-random Borel sets $E\subseteq[0,1]$ that 
intersect $\Lambda_H$. Our characterization 
is not so easy to describe here in the Introduction. 
For now, suffice it to say that 
it readily yields Theorem~\ref{thm:Hdim}. 
The following is another 
consequence of the said characterization: If 
$\nu$ is standard normal, then 
\begin{equation}\label{eq:Pdim} 
   \sup_{t\in E} \limsup_{n\to\infty} 
   \frac{\left( S_n(t) \right)^2 - 2n\ln\ln n}{n\ln\ln\ln n} 
   = 3 + 2\Pdim E. 
\end{equation} 
Here, $\Pdim$ denotes packing dimension~\cite{Mattila}. 
The preceding display follows from \eqref{eq:Pdim1} below. 
 
On one hand, if we set $E$ to be the 
entire interval $[0,1]$, then the right-hand side of 
\eqref{eq:Pdim} is 
equal to $5$, and we obtain an earlier result of 
ours~\ycite{KLM}*{Eq.~1.15}. On the other hand, 
if we set $E$ to be a singleton, then the right-hand side 
of (\ref{eq:Pdim}) 
is equal to $3$, and we obtain the second-term correction 
to the classical law of the iterated 
logarithm~\citelist{\cite{kolmogorov}\cite{erdos}}. 
 
Somewhat unexpectedly, the next result follows also 
from Theorem~\ref{thm:genest}. To the best of our 
knowledge it is new. 
 
\begin{corollary}\label{cor:genest-OU} 
   Let $\{Z_t\}_{t\ge 0}$ denote the Ornstein--Uhlenbeck (OU) process 
   on the real line that satisfies the s.d.e.\@ 
   $dZ = -Z\, dt +\sqrt{2}\, dW$ for a Brownian motion $W$. 
   Then for every non-empty non-random closed set $E\subseteq[0,1]$, 
   and all $z>1$, 
   \begin{equation} 
      A_{\ref{eq:genest}}^{-1} \K_E\left( 
      \frac{1}{z^2}\right) \bar\Phi(z) \le 
      \P\left\{ \sup_{t\in E} Z_t \ge z 
      \right\} \le A_{\ref{eq:genest}} 
      \K_E\left( \frac{1}{z^2}\right) \bar\Phi(z). 
   \end{equation} 
\end{corollary} 
Section~\ref{sec:OU} below contains further remarks along these lines. 
 
For a proof of Corollary~\ref{cor:genest-OU} consider 
the two-parameter processes, 
\begin{equation} 
   \mathscr{S}_n(u,t) = \frac{X_1(t)+\cdots + X_{\lfloor 
   un\rfloor} (t)}{\sqrt n}\quad (0\le u,t\le 1),\ 
   {}^\forall n=1,2,\ldots\,. 
\end{equation} 
Our recent work~\ycite{KLM}*{Theorem 1.1} 
implies that if $\nu$ is standard 
normal, then $\mathscr{S}_n\Rightarrow \mathscr{U}$ 
in the sense of $\mathscr{D}([0,1]^2)$~\cite{bickel}, and 
$\mathscr{U}$ is the continuous centered Gaussian process 
with correlation function 
\begin{equation}\label{eq:covarianceU} 
   \E\left[ \mathscr{U}(u,t) 
   \mathscr{U}(v,s) \right] = e^{-|t-s|} \min(u,v), 
   \quad 0\le u,v,s,t\le 1. 
\end{equation} 
In particular, $\sup_{t\in E} S_n(t)/\sqrt n$ converges in 
distribution to $\sup_{t\in E} \mathscr{U}(1,t)$. 
Corollary~\ref{cor:genest-OU} follows from this and the fact that 
$\sup_{t\in E} \mathscr{U}(1,t)$ has the same distribution 
as $\sup_{t\in E} Z_t$. 
 
In this paper we apply stochastic calculus to strengthen our 
earlier central limit theorem~\ycite{KLM}*{Theorem 1.1}. 
Indeed we offer the following invariance principle. 
 
\begin{theorem}\label{thm:invariance} 
   If $\int_{-\infty}^\infty x\, \nu(dx)=0$
   and $\int_{-\infty}^\infty x^2\, \nu(dx)=1$,
   then $\mathscr{S}_n\Rightarrow \mathscr{U}$ in the 
   sense of $\mathscr{D}([0,1]^2)$. 
\end{theorem} 
 
We close the introduction by presenting the following 
dynamic stability result. 
 
\begin{theorem}\label{thm:DS-REC} 
   Suppose $\nu$ is a distribution on $\mathbf{Z}$ 
   which has mean zero and variance one. If there exists 
   $\e>0$ such that $\int_{-\infty}^\infty 
   |x|^{2+\e}\, \nu(dx)<\infty$, then 
   \begin{equation} 
      \P\left\{ \sum_{n=1}^\infty \mathbf{1}_{ 
      \{S_n(t)=0\}} =\infty\ \text{ for all } 
      t\ge 0\right\} =1. 
   \end{equation} 
\end{theorem} 
 
In words, under the conditions of Theorem~\ref{thm:DS-REC}, 
the recurrence of the origin 
is dynamically stable. 
When $\nu$ is supported by a finite subset of $\mathbf{Z}$ 
this was proved by~\ocite{benjamini}*{Theorem 1.11}. In 
order to generalize to the present setting, we first 
develop the following quantitative form 
of the classical gambler's ruin theorem. We state it 
next, since it may be 
of independent interest. 
 
Consider i.i.d.\@ integer-valued 
random variables $\{\xi_n\}_{n=1}^\infty$ such that 
$\E[\xi_1]=0$ and $\s^2=\E[\xi_1^2]<\infty$. Define 
$s_n=\xi_1+\cdots+\xi_n$ to be the corresponding random 
walk, and let $T(z)$ denote the first-passage time to $z$; i.e., 
\begin{equation} 
   T(z)= \inf\left\{ n\ge 1: s_n =z\right\} 
   \qquad {}^\forall z\in\mathbf{Z} 
   \quad(\inf\varnothing=\infty). 
\end{equation} 
 
\begin{theorem}[Gambler's Ruin]\label{thm:gamblersruin} 
   If $G$ denotes the additive subgroup 
   of $\mathbf{Z}$ generated by the possible values of 
   $\{s_n\}_{n=1}^\infty$, then there 
   exists a constant $A_{\ref{eq:gamblersruin}} 
   =A_{\ref{eq:gamblersruin}}(\s,G)>1$ such that 
   \begin{equation}\label{eq:gamblersruin} 
      \frac{A_{\ref{eq:gamblersruin}}^{-1}}{1+|z|} 
      \le \P\left\{ T(z) \le T(0) \right\} \le 
      \frac{A_{\ref{eq:gamblersruin}}}{1+|z|} 
      \qquad{}^\forall z\in G. 
   \end{equation} 
\end{theorem} 
\bigskip 
\noindent\textbf{Acknowledgements}\ We wish to thank 
Professor Harry Kesten for discussions regarding
Theorem~\ref{thm:gamblersruin}, and Professor
Mikhael Lifshits for bringing the work of
\ocite{rusakov:95} to our attention.
 
\section{On the Kolmogorov $\e$-Entropy} 
 
\subsection{$\Lambda_H$-Polar Sets} 
Let $H:\R_+\to\R_+$ be non-decreasing and measurable, 
and recall the random set $\Lambda_H$ from (\ref{eq:KLM}). 
 
We say that a measurable set $E\subset[0,1]$ is 
\emph{$\Lambda_H$-polar} if $\P\{\Lambda_H\cap E\neq\varnothing\}=0$. 
If $E$ is not $\Lambda_H$-polar, then the Hewitt--Savage law insures 
that $\P\{\Lambda_H\cap E\neq\varnothing\}=1$. Our characterization 
of $\Lambda_H$-polar sets is described in terms of the function 
\begin{equation}\label{eq:psi} 
   \psi_H(E) = 
   \int_1^\infty H^2(t) \K_E \left( \frac{1}{H^2(t)} 
   \right) \bar\Phi (H(t))\, \frac{dt}{t},\qquad 
   {}^\forall E\subseteq[0,1]. 
\end{equation} 
Although $\psi_H$ is subadditive, 
it is not a measure; e.g., $\psi_H$ assigns equal 
mass $\int_1^\infty H^2(t)\bar\Phi(H(t))\,\frac{dt}{t}$ 
to all singletons. We will show that the function 
$\psi_E$ determines the growth-rate 
of $\sup_{t\in E}S_n(t)$ in the following sense. 
 
\begin{theorem}\label{thm:inttest1} 
   Suppose $E\subseteq [0,1]$ is Borel-measurable 
   and $H:\R_+\to\R_+$ is non-decreasing. Then, 
   \begin{equation} 
      \limsup_{n\to\infty} 
      \left[ \sup_{t\in E} S_n(t) - H(n)\sqrt{n} \right] > 0\ 
      \text{if and only if $\psi_H(E)=\infty$}. 
   \end{equation} 
\end{theorem}

\begin{remark}\label{rem:inttest1} 
   In fact, we will prove that: 
   \begin{equation}\begin{split} 
      \psi_H(E) =\infty & \Longrightarrow 
          \limsup_{n\to\infty} 
          \left[ \sup_{t\in E} S_n(t) - H(n)\sqrt{n} \right] 
          =\infty;\\ 
      \psi_H(E)<\infty & \Longrightarrow 
          \limsup_{n\to\infty} 
          \left[ \sup_{t\in E} S_n(t) - H(n)\sqrt{n} \right] 
          =-\infty. 
   \end{split}\end{equation} 
\end{remark} 
 
\begin{definition} 
   We write $\Psi_H(E)<\infty$ if we 
   can decompose $E$ as $E=\cup_{n=1}^\infty E_n$---where 
   $E_1,E_2,\ldots,$ are closed---such that for all $n\ge 1$, 
   $\psi_H(E_n)<\infty$. Else, we say that $\Psi_H(E)=\infty$. 
\end{definition} 
\begin{remark}\label{rem:psiPsi} 
   One can have $\Psi_H(E)<\infty$ although 
   $\psi_H(E)=\infty$. See Example~\ref{ex:psiPsi} below. 
\end{remark} 
 
The following then characterizes all polar sets of 
$\Lambda_H$; it will be shown to 
be a ready consequence of Theorem~\ref{thm:inttest1}. 
 
\begin{theorem}\label{thm:inttest2} 
   Suppose $E\subset [0,1]$ is a fixed compact set, 
   and $H:\R_+\to\R_+$ is non-decreasing. Then, 
   $E$ is $\Lambda_H$-polar if and only if $\Psi_H(E)= \infty$. 
\end{theorem} 
 
\begin{remark}\label{rem:inttest2} 
   The following variation of Remark~\ref{rem:inttest1} 
   is valid: 
   \begin{equation}\begin{split} 
      \Psi_H(E) =\infty & \Longrightarrow 
          \sup_{t\in E} \limsup_{n\to\infty} 
          \left[ S_n(t) - H(n)\sqrt{n} \right] 
          =\infty;\\ 
      \Psi_H(E)<\infty & \Longrightarrow 
          \sup_{t\in E} \limsup_{n\to\infty} 
          \left[ S_n(t) - H(n)\sqrt{n} \right] 
          =-\infty. 
   \end{split}\end{equation} 
\end{remark} 
 
\subsection{Relation to Minkowski Contents} 
In the remainder of 
this section we say a few words about the function 
$\K_E$. To begin with, let us note that the defining 
\emph{maximal Kolmogorov sequence} $\{x_j\}_{j=1}^k$ 
has the property that any 
\begin{equation} \label{kolmo} 
  w\in E  \text{ satisfies }  |w-x_j|\le\e \text{ for some } j=1,\ldots,k. 
\end{equation}

The Kolmogorov $\e$-entropy 
is related to the \emph{Minkowski content} 
of $E\subseteq \R$. The latter can be defined as follows: 
\begin{equation} 
   \M_n(E) = \sum_{i=-\infty}^\infty a_{i,n}(E), 
   \text{ where }a_{i,n}(E) = 
   \begin{cases} 
      1, & \text{ if $\left[ \frac{i}{n},\frac{i+1}{n}\right) 
         \cap E\neq\varnothing$},\\ 
      0, &\text{ otherwise}. 
   \end{cases} 
\end{equation} 
 
Here is the relation. 
See~\ocite{dudley}*{Theorem 6.0.1} and~\ocite{Mattila}*{p.\@ 78, eq.\@ 5.8} 
for a related inequality. 
 
\begin{proposition}\label{pr:KM} 
   For all non-empty sets $E\subseteq[0,1]$ and all integers $n\ge 1$, 
   \begin{equation} 
      \K_E (1/n) \le \M_n(E) \le 3\K_E (1/n). 
   \end{equation} 
\end{proposition} 
 
\begin{remark}
   It is not difficult to see that
   both bounds can be attained.
\end{remark}

\begin{proof} 
   Let $k=\K_E(1/n)$ and choose  maximal (Kolmogorov) 
   points $x_1 < x_2\ldots <x_k$ such that any distinct pair $(x_i,x_j)$ are 
   distance at least $1/n$ apart. Define $\mathscr{E}$
   to be the collection of all intervals $[\frac{i}{n},\frac{i+1}{n})$,
   $0\le i<n$, such that any $I\in\mathscr{E}$ intersects $E$.
   Let $\mathscr{G}$ denote the collection of all $I\in\mathscr{E}$
   such that some $x_j$ is in $I$. These are the ``good'' intervals.
   Let $\mathscr{B} = \mathscr{E}\setminus\mathscr{G}$ denote
   the ``bad'' ones. Good intervals contain exactly
   one of the maximal Kolmogorov points, whereas bad ones contain 
   none. 
   Therefore, $\K_E(1/n)= |\mathscr{G}| \le |\mathscr{E}| =\M_E(1/n)$,
   where $|\cdots|$ denotes cardinality.
   To complete our derivation we prove that $|\mathscr{B}|\le 2\K_E(1/n)$.
   
   We observe that any bad interval is necessarily adjacent to
   a good one. Therefore, we can
   write $\mathscr{B}=\mathscr{B}_L \cup \mathscr{B}_R$ where
   $\mathscr{B}_L$ [resp.\@ $\mathscr{B}_R$]
   denotes the collection of all bad intervals
   $I$ such that there exists a good interval adjacent to the left
   [resp.\@ right] of $I$. By virtue of their definition,
   both $\mathscr{B}_L$ and $\mathscr{B}_R$ each have 
   no more than $\K_E(1/n)$ elements. This completes
   the proof.
\end{proof} 
   An immediate consequence of this result is that if $\e \in [ 
    2^{-n-1}, 2^{-n}]$ then 
 \begin{equation}\label{eq:K=RV} 
      \K_E(\e) \le \K_E\left( 2^{-n-1} \right) 
      \le M_{2^{n+1}}(E) \leq 2 M_{2^{n}}(E)
      \leq 6  \K_E\left( 2^{-n} \right).
 \end{equation} 
 
\subsection{Relation to Minkowski and Packing Dimensions} 
 
There are well-known connections between $\e$-entropy and 
the (upper) Minkowski dimension, some of which 
we have already seen; many more of which one can find, 
in fine pedagogic form, in~\ocite{Mattila}*{Ch.\@ 5}. 
We now present a relation that is particularly 
suited for our needs. Let $H_\rho$ be any locally-bounded 
non-decreasing function such that 
\begin{equation} 
   H_\rho(t) = \sqrt{ 2\ln \ln t + 
   2\rho \ln \ln \ln t},\qquad{}^\forall t>e^{10000}. 
\end{equation} 
One or two lines of calculations then reveal that 
\begin{equation}\label{eq:asymp} 
   \psi_{H_\rho}(E) <\infty\ \text{if and only if}\ 
   \int_1^\infty \K_E(1/s) s^{\frac12 -\rho}\, ds <\infty. 
\end{equation} 
 
\begin{proposition}\label{pr:EM} 
   For all compact linear sets $E$, 
   \begin{equation}\begin{split} 
      \Kdim E &= \inf\left\{\rho>0:\ 
         \psi_{H_\rho}(E) < \infty\right\} -\frac32,\text{ and}\\ 
      \Pdim E  &= \inf\left\{\rho>0:\ 
         \Psi_{H_\rho}(E) < \infty\right\} -\frac32. 
   \end{split}\end{equation} 
\end{proposition} 
 
There are well-known examples of sets $E$ 
whose packing and upper Minkowski dimension 
differ. Therefore, Proposition~\ref{pr:EM} provides 
us with an example of functions $H$ (namely an 
appropriate $H_\rho$) and sets $E$ such that 
$\psi_H(E)$ is infinite although 
$\Psi_H(E)$ is finite. This is good enough to 
address the issue raised in Remark~\ref{rem:psiPsi}. 
In fact, one can do more at little extra cost. 
 
\begin{example}\label{ex:psiPsi} 
   Define 
   \begin{equation}\label{eq:zeta} 
      \mathscr{J}_\zeta (H) = \int_1^\infty 
      H^\zeta(t)\bar\Phi(H(t))\,\frac{dt}{t}\qquad 
      {}^\forall \zeta>0. 
   \end{equation} 
   Now consider any measurable non-decreasing 
   function $H:\R_+\to\R_+$ such that 
   $\mathscr{J}_2(H)<\infty$ but $\mathscr{J}_{2+\e}(H)=\infty$ 
   for some $\e>0$. 
   Then there are compact sets $E\subseteq[0,1]$ 
   such that $\psi_H(E)=\infty$ although 
   $\Psi_H(E)<\infty$. Our construction of such 
   an $E$ is based on a well-known example 
   \cite{Mattila}*{Exercise 1, p.\@ 88}. 
 
   Without loss of generality, we may assume 
   that $\e\in(0,1)$. Bearing this in mind, define $r_0=1$ and 
   $r_k = 1 - \sum_{j=1}^k j^{-1/\e}$ ($k=1,2,\ldots$). 
   Now consider 
   \begin{equation} 
      E = \{0\} \cup \bigcup_{k=0}^\infty \{ r_k\}. 
   \end{equation} 
   Then it is possible to prove that there is a constant 
   $A>1$ such that for all $\d\in(0,1)$, 
   $A^{-1} \d^\e \le \K_E(\d) \le A \d^\e$. 
   In particular, $\psi_H(E)$ is comparable to 
   $\mathscr{J}_{2+\e}(H)=\infty$. On the other hand, 
   because $E$ is countable and $\mathscr{J}_2(H)<\infty$, 
   we readily have $\Psi_H(E)<\infty$. 
\end{example} 
 
Our proof of Proposition~\ref{pr:EM} requires the following 
little lemma from geometric measure theory. 
 
\begin{lemma}\label{lem:Psi} 
   Suppose $H:\R_+\to\R_+$ is non-decreasing and measurable, 
   and $E\subseteq[0,1]$ is Borel and satisfies 
   $\Psi_H(E)=\infty$. Then, there exists a compact 
   set $G\subseteq E$ such that $\psi_H(I\cap G)=\infty$ 
   for all rational intervals $I\subseteq[0,1]$ that intersect $G$. 
\end{lemma} 
 
\begin{proof} 
   Let $\mathscr{R}$ denote the collection of all 
   open rational intervals in $[0,1]$, and define 
   \begin{equation} 
      E_*=\bigcup_{I\in\mathscr{R}:~\psi_H(E\cap I)<\infty}I. 
   \end{equation} 
   A little thought makes it manifest that 
   $E_*$ is an open set in $[0,1]$, and 
   $G=E\setminus E_*$ has the desired properties. 
\end{proof} 
 
\begin{proof}[Proof of Proposition~\ref{pr:EM}] 
   We will prove the assertion about $\Kdim$; the 
   formula for $\Pdim$ follows from the one 
   for $\Kdim$, Lemma~\ref{lem:Psi}, 
   and regularization~\cite{Mattila}*{p.\@ 81}. 
 
   Throughout the proof, we let 
   $d=\Kdim(E)$ denote the Minkowski dimension of 
   $E$~\cite{Mattila}*{p.\@ 79}. 
   By its very definition, and thanks to 
   Proposition~\ref{pr:KM}, $d$ can be written as 
   \begin{equation} 
      d = \Kdim E 
      = \limsup_{s\to\infty} \frac{\log \K_E(1/s)}{\log s}. 
   \end{equation} 
 
   Now 
   \begin{equation}\begin{split} 
      \int_1^\infty \K_E(1/s) s^{\frac12 -\rho}\, ds 
         & = \sum_{n=0}^\infty \int_{2^n}^{2^{n+1}} 
         \K_E(1/s) s^{\frac12 -\rho}\, ds\\ 
      & \ge 2^{-\rho} 
         \sum_{n=0}^\infty \K_E\left( 2^{-n}\right) 
         2^{-(\rho-\frac32)n}\\ 
      & \ge 2^{-\rho}\limsup_{n\to\infty} \K_E\left( 2^{-n}\right) 
         2^{-(\rho-\frac32)n}. 
   \end{split}\end{equation} 
   Thus, if $2^n\le s\le 2^{n+1}$ and $\rho>2$, then 
   for all sufficiently large $n$, 
   \begin{equation} 
      s^{-(\rho-\frac32)} \K_E(1/s) \le 6\cdot 
      2^{-n(\rho-\frac32)} \K_E\left( 2^{-n} \right). 
   \end{equation} 
   See (\ref{eq:K=RV}). This development shows that 
   \begin{equation} 
      \int_1^\infty \K_E(1/s) s^{\frac12-\rho}\, ds 
      \ge \frac{1}{6 \cdot 2^\rho} \limsup_{s\to\infty} 
      \frac{\K_E(1/s)}{s^{(\rho-\frac32)}}. 
   \end{equation} 
   Therefore, whenever $\rho-\frac32<d$, the integral 
   on the left-hand side is infinite. Thanks to 
   (\ref{eq:asymp}), this means that 
   \begin{equation} 
      \inf\left\{ \rho>0:\ \psi_{H_\rho}(E) 
      <\infty\right\} \le \frac32+d = \frac32+\Kdim E. 
   \end{equation} 
   This is half of the result for the Minkowski dimension. 
   To prove the converse half, we argue similarly, and 
   appeal to (\ref{eq:K=RV}), to deduce that 
   \begin{equation} 
      \int_1^\infty \K_E(1/s) s^{\frac12-\rho}\, ds  \le 
      6 \sum_{n=0}^\infty \K_E\left( 2^{-n} \right) 
      2^{-n(\rho-\frac32)} \le 
      6 \sum_{n=0}^\infty 2^{n(d-\rho+\frac32) +o(n)}. 
   \end{equation} 
   In particular, if $\rho>d+\frac32$, then the left-hand 
   side is finite. This and (\ref{eq:asymp}) 
   together verify the asserted identity for $\Kdim$. 
\end{proof} 
 
\begin{remark} 
   In conjunction, 
   Theorem~\ref{thm:inttest2} and Proposition~\ref{pr:EM} 
   show that for any non-random Borel set $E\subseteq[0,1]$, 
   \begin{equation}\label{eq:pdimchar}\begin{split} 
      \rho> \frac32 +\Pdim E &\ \Longrightarrow\ 
         \Lambda_{H_\rho} \cap E = \varnothing\\ 
      \rho< \frac32 +\Pdim E &\ \Longrightarrow\ 
         \Lambda_{H_\rho} \cap E \neq \varnothing. 
   \end{split}\end{equation} 
   Moreover, the intersection argument of~\ocite{KPX}*{Theorem 3.2} 
   goes through unhindered to imply that if 
   $\rho < \frac32+\Pdim E$, then $\Pdim(\Lambda_H\cap E)=\Pdim E$. 
   In particular, we can apply this with $E=[0,1]$, and recall 
   (\ref{KLM:T1.5}), to deduce the following: 
   \begin{equation}\label{eq:Pdim1}\begin{split} 
      \rho &< \frac52\ \Longrightarrow\ \Pdim \Lambda_{H_\rho} =1,\\ 
      \rho &> \frac52\ \Longrightarrow\ \Lambda_{H_\rho}=\varnothing. 
   \end{split}\end{equation} 
   Equation (\ref{eq:Pdim}) is an immediate consequence of this. 
   One could alternatively use the 
   limsup-random-fractal theories of~\ocite{KPX} and~\ocite{DPRZ} 
   to derive (\ref{eq:Pdim1}). 
\end{remark} 
\subsection{An Application to Stable Processes}\label{subsec:stable} 
 
Let $\{Y_\alpha(t)\}_{t\ge 0}$ denote a symmetric 
stable process with index 
$\alpha\in(0,1)$, and let us consider the random set 
$\mathscr{R}_\alpha =\text{cl} (Y_\alpha([1,2]))$ 
denote the closed range of $\{Y_\alpha(t)\}_{t\in[1,2]}$.

\begin{proposition}\label{pr:stableindex} 
   Consider a given $\alpha,\beta\in(0,1)$. Then, 
   for all $M>0$ and $p\ge 1$, 
   there exists a finite constant 
   $A_{\ref{eq:stableindex}}=A_{\ref{eq:stableindex}} 
   (\alpha,\beta,p,M)>1$ such that for all intervals $I\subset[-M,M]$ 
   with length $\ge\beta$, and all $\e\in(0,1)$, 
   \begin{equation}\label{eq:stableindex} 
      A_{\ref{eq:stableindex}}^{-1}\e^{-\alpha p} 
      \le \E\left[ \K^p_{\mathscr{R}_\alpha\cap I} (\e) 
      \right]\le 
      A_{\ref{eq:stableindex}}\e^{-\alpha p}. 
   \end{equation} 
\end{proposition} 
 
\begin{proof} 
   Thanks to Proposition~\ref{pr:KM}, it suffices 
   to show that we can find $A_{\ref{eq:stableindex2}}>1$ 
   [depending only on $\alpha,M,p$] 
   such that for all $n\ge 1$, 
   \begin{equation}\label{eq:stableindex2} 
      A_{\ref{eq:stableindex2}}^{-1} n^{\alpha p} 
      \le \E\left[ \M^p_n \left( \mathscr{R}_\alpha\cap I\right) 
      \right]\le 
      A_{\ref{eq:stableindex2}}n^{\alpha p}. 
   \end{equation} 
 
   This follows from connections 
   to potential-theoretic notions, for which 
   we need to introduce some notation. 
 
   Let $p_t(x,y)$ denote the transition densities of 
   the process $Y_\alpha$. As usual, $P_x$ denotes the law 
   of $x+Y_\alpha(\bullet)$ on path-space. 
   Define $r(x,y)$ to be the $1$-potential density 
   of $Y_\alpha$; i.e., 
   \begin{equation} 
      r(x,y) = \int_0^\infty e^{-s} p_s(x,y)\, ds. 
   \end{equation} 
   Finally, let $T(z,\e) = \inf \{ s>0:\ 
   |Y_\alpha(s) - z|\le \e \}$ designate the entrance 
   time of the interval $[z-\e,z+\e]$; as usual, 
   $\inf\varnothing =\infty$. 
 
   It is well known that for any $M>0$, 
   there exists a constant 
   $A=A(M,\alpha)>1$ such that 
   \begin{equation}\begin{split} 
      A^{-1} \e^{1-\alpha} &\le 
         \inf_{x\in[-M,M]} 
         \P \left\{ \mathscr{R}_\alpha\cap [x-\e,x+\e] 
         \neq \varnothing \right\}\\ 
      &  \le \sup_{x\in\R} 
         \P \left\{ \mathscr{R}_\alpha\cap [x-\e,x+\e] 
         \neq \varnothing \right\} \le A\e^{1-\alpha}; 
   \end{split}\end{equation} 
   see, for example \ocite{Khoshnevisan}*{Proposition 1.4.1, p.\@ 351}. 
   In the case $p=1$, this proves 
   Equation~(\ref{eq:stableindex2}). 
   Because $L^p(\P)$-norms are increasing in $p$, 
   the lower bound in (\ref{eq:stableindex2}) follows, 
   in fact, for all $p\ge 1$. Thus, it remains to prove 
   the corresponding upper bound for $p>1$. 
 
   Modern variants of classical probabilistic potential 
   theory tell us that for all $x\not\in[y-\e,y+\e]$, 
   \begin{equation}\begin{split} 
      &\int_0^\infty e^{-s} P_x \left\{ 
         T(y,\e) \le s\right\} \, ds\\ 
      &\quad \le \mathbf{S} 
         \left[ \inf_{\mu\in\mathscr{P}([y-\e,y+\e])} 
         \iint r(u,v)\, \mu(du)\,  \mu(dv)\right]^{-1}. 
   \end{split}\end{equation} 
   See~\ocite{Khoshnevisan}*{Theorem 2.3.1, P.\@ 368}. 
   Here, $\mathbf{S} = \sup_{z\in  [y-\e,y+\e]} r(x,z)$, 
   In the preceding, $E$ is a linear Borel 
   set, and $\mathscr{P}(E)$ denotes the collection 
   of all probability measures on the Borel set $E$. 
 
   On the other hand, there exists a finite constant 
   $A_{\ref{eq:r-est}}>1$ such that whenever $x,y$ are both in $[-2M,2M]$, 
   \begin{equation}\label{eq:r-est} 
     A_{\ref{eq:r-est}}^{-1} |x-y|^{-1+\alpha} \le 
     r(x,y) \le A_{\ref{eq:r-est}} |x-y|^{-1+\alpha}. 
   \end{equation} 
   See, for example,~\ocite{Khoshnevisan}*{Lemma 3.4.1, p.\@ 383}. 
   Now as soon as we have $|x-y|\ge 2\e$ and $|z-y|\le\e$, 
   it follows that $|x-z|\ge\frac12 |x-y|$. 
   Therefore, the inequality $\int_0^\infty (\cdots) 
   \ge \int_0^1 (\cdots)$ leads us to the following: 
   \begin{equation}\begin{split} 
     &P_x \left\{  T(y,\e) \le 1\right\}\\ 
     &\quad \le 2^{1-\alpha} eA_{\ref{eq:r-est}}^2 |x-y|^{-1+\alpha} 
        \left[ \inf_{\mu\in\mathscr{P}([-\e,+\e])} 
        \iint |u-v|^{-1+\alpha} 
        \, \mu(du)\,  \mu(dv)\right]^{-1}. 
   \end{split}\end{equation} 
   The term $[\cdots]^{-1}$ is the $(1-\alpha)$-dimensional 
   Riesz capacity of $[-\e,\e]$. It is a classical fact 
   that the said capacity is, up to multiplicative constants, 
   of exact order $\e^{1-\alpha}$. 
   Therefore, there exists $A_{\ref{eq:hit-UB}}>1$ 
   such that for all $\e\in(0,1)$ 
   and all $x,y\in[-2M,2M]$ that satisfy 
   $|x-y|\ge 2\e$, 
   \begin{equation}\label{eq:hit-UB} 
      P_x \left\{  T(y,\e) \le 1\right\} 
     \le A_{\ref{eq:hit-UB}} |x-y|^{-1+\alpha} \e^{1-\alpha}. 
   \end{equation} 
   We now prove the upper bound in (\ref{eq:stableindex2}) 
   for the case $p=2$ and hence all $p\in[1,2]$. 
   By the strong Markov property and time reversal, 
   whenever $x,y\in[-2M,2M]$ satisfy $|x-y|\ge 4\e$, 
   \begin{equation}\begin{split} 
      &\P\left\{ \mathscr{R}_\alpha \cap [x-\e,x+\e] \neq 
         \varnothing ~,~ \mathscr{R}_\alpha \cap [y-\e,y+\e] 
         \neq \varnothing \right\}\\ 
      &\quad \le 2 \P\left\{ \mathscr{R}_\alpha \cap [x-\e,x+\e] \neq 
         \varnothing \right\} \sup_{v\in[x-\e,x+\e]} P_v \left\{ 
         T(y,\e) \le 1 \right\}\\ 
      &\quad \le 2A_{\ref{eq:hit-UB}} |x-y|^{-1+\alpha}\e^{2(1-\alpha)}. 
   \end{split}\end{equation} 
   Equation~(\ref{eq:stableindex}) readily follows from this in the case that 
   $p=2$. To derive the result for an arbitrary positive integer 
   $p$, simply iterate this argument $p-1$ times. 
\end{proof} 
\section{Proof of Theorem~\ref{thm:genest}} 
 
This proof rests on half of the following 
preliminary technical result. 
Throughout this section $\{z_n\}_{n=1}^\infty$ 
is a fixed sequence that satisfies (\ref{eq:z}), 
and $E \subseteq [0,1]$ is a fixed non-random compact set. 
 
\begin{proposition}\label{pr:TechnicalEst} 
   Let $\{\d_n\}_{n=1}^\infty$ be a fixed sequence of numbers in 
   $[0,1]$ that satisfy 
   \begin{equation}\label{eq:DeltaGap} 
      \liminf_{n\to\infty}\d_n z_n^2>0. 
   \end{equation} 
   Then there exists a finite constant $A_{\ref{eq:TechnicalEst}}>1$ 
   such that for all $n\ge 1$, 
   \begin{equation}\label{eq:TechnicalEst} 
      A_{\ref{eq:TechnicalEst}}^{-1} 
      \d_n z_n^2\bar\Phi(z_n) \le 
      \P\left\{ \sup_{t\in[0,\d_n]} S_n(t) \ge z_n\sqrt{n} \right\} 
      \le A_{\ref{eq:TechnicalEst}} \d_n z_n^2\bar\Phi(z_n). 
   \end{equation} 
\end{proposition} 
 
\begin{proof} 
   We will need some of the notation, 
   as well as results, of~\ocite{KLM}. Therefore, we 
   first recall the things that we need. 
 
   Let $\PN$ (resp.\@ $\EN$) denote the `quenched' 
   measure $\P(\cdots\,|\,\mathscr{N})$ 
   (resp.\@ expectation operator $\E[\cdots\,|\, \mathscr{N}]$), where 
   $\mathscr{N}$ denotes the $\s$-algebra generated by all 
   of the clocks, and define $\F^n_t$ to be the $\s$-algebra 
   generated by $\{S_j(s);\ 0\le s\le t\}_{j=1}^n$. 
 
   Define 
   \begin{equation}\begin{split} 
      L_n(t) &= \int_0^t \mathbf{1}_{B_n(u) }\, du,\text{ where}\\ 
      B_n(t) &= \left\{ \omega\in\Omega:\ 
         S_n(t) \ge z_n\sqrt{n} \right\}. 
   \end{split}\end{equation} 
   We replace the variable $J_n$ of~\ocite{KLM}*{eq.\@ 5.3} 
   by our $L_n(2\d_n)$, and go through the proof of~\ocite{KLM}*{Lemma 5.2} 
   to see that there exists 
   an $\mathscr{N}$-measurable event $A_{n,\frac12}$ such that 
   for any $u\in[0,\d_n]$, the following holds 
   $\P$-almost surely: 
   \begin{equation}\label{eq:CondBd}\begin{split} 
      \EN\left[ \left. L_n(2\d_n) \, \right|\, \F^n_u \right] 
         & \ge \frac{2}{3 z_n^2} \int_0^{\frac32 (2\d_n-u) z_n^2} 
         \bar\Phi \left(\sqrt{t}\right)\, dt \cdot 
         \mathbf{1}_{A_{n,1/2} \cap B_n(u)}\\ 
      & \ge \frac{2}{3 z_n^2} \int_0^{\frac32 \d_n z_n^2} 
         \bar\Phi \left(\sqrt{t}\right)\, dt \cdot 
         \mathbf{1}_{A_{n,1/2} \cap B_n(u)}\\ 
      & \ge \frac{A_{\ref{eq:CondBd}}}{z_n^2}\cdot 
         \mathbf{1}_{A_{n,1/2} \cap B_n(u)}, 
   \end{split}\end{equation} 
   where $A_{\ref{eq:CondBd}}$ is an absolute constant that is 
   bounded below. 
   Moreover, thanks to~\ocite{KLM}*{Theorem 2.1} and (\ref{eq:DeltaGap}), there exists 
   a finite constant $A_{\ref{eq:AC}}\in(0,1)$ such that for all $n\ge 1$, 
   \begin{equation}\label{eq:AC} 
      \P\left( A_{n,\frac12}^\complement \right) 
      \le  \,z_n^2 \delta_n\, e^{-A_{\ref{eq:AC}} 
      n/z_n^2}. 
   \end{equation} 
   Now, $u\mapsto \EN[L_n(2\d_n)\,|\,\F^n_u]$ is a non-negative 
   and bounded $\PN$-martingale. Therefore, $\P$-almost surely, 
   \begin{align} 
      \mathbf{1}_{A_{n,1/2}} \PN \left\{ 
         \sup_{u\in[0,\d_n]} S_n(t) \ge z_n\sqrt{n}\right\} 
         & = \PN \left\{ \sup_{u\in[0,\d_n]\cap\mathbf{Q} } 
         \mathbf{1}_{ A_{n,1/2} \cap B_n(u)} \ge 1 
         \right\}\nonumber\\ 
      & \le \PN\left\{ \sup_{u\in[0,\d_n]\cap\mathbf{Q} } 
         \EN\left[ L_n(2\d_n)\, \big|\, \F^n_u \right] 
         \ge \frac{A_{\ref{eq:CondBd}}}{z_n^2} \right\}\\ 
      & \le \frac{z_n^2}{A_{\ref{eq:CondBd}}} 
         \EN\left[ L_n(2\d_n)\right] 
         = \frac{2}{A_{\ref{eq:CondBd}}} 
         \d_n z_n^2\bar\Phi(z_n).\nonumber 
   \end{align} 
   The ultimate inequality follows from Doob's maximal inequality 
   for martingales, 
   and the last equality from the stationarity of $t\mapsto S_n(t)$. 
   Taking expectations and applying (\ref{eq:AC}) yields 
   \begin{equation} 
      \P\left\{ \sup_{t\in[0,\d_n]} S_n(t) 
      \ge z_n\sqrt{n} \right\} \le 
      \frac{2}{A_{\ref{eq:CondBd}}} 
      \d_n z_n^2 \left[\,\bar\Phi(z_n) + e^{-A_{\ref{eq:AC}} 
      n/z_n^2}.\,\right] 
   \end{equation} 
   Equation (\ref{eq:z}) shows that the first term 
   on the right-hand side dominates the second one for all $n$ 
   sufficiently large. This yields the probability upper bound of 
   the proposition. Now we work toward the lower bound. 
 
   By adapting the argument of~\ocite{KLM}*{eq.\@ 6.12}, we can conclude 
   that $\P$-almost surely there exists an $\mathscr{N}$-measurable 
   $\P$-a.s. finite random variable $\gamma$ such that for all $n\ge \gamma$, 
   \begin{equation}\label{eq:L2L} 
      \EN\left[ \left( L_n(\d_n) \right)^2 \right] \le 
      A_{\ref{eq:L2L}} \d_n z_n^{-2} \bar\Phi(z_n). 
   \end{equation} 
   where $A_{\ref{eq:L2L}}>1$ is a non-random and finite 
   constant. [Replace $J_n$ by $L_n(\d_n)$ and 
   proceed to revise equation (6.12) of~\ocite{KLM}.] 
   Since, by stationarity, 
   $\EN[L_n(\d_n)]=\d_n\bar\Phi(z_n)$, 
   the Paley--Zygmund inequality shows that $\P$-almost surely 
   for all $n\ge \gamma$, 
   \begin{equation} 
      \PN\left\{ L_n(\d_n)>0\right\} \ge 
      \frac{\left(\EN[L_n(\d_n)]\right)^2}{\EN\left[ 
      \left( L_n(\d_n) \right)^2\right]} \ge 
      \frac{1}{A_{\ref{eq:L2L}}} \d_n z_n^2 \bar\Phi(z_n). 
   \end{equation} 
   On the other hand, 
   \begin{equation} 
      \P\left\{ \sup_{t\in[0,\d_n]} S_n(t)\ge 
      z_n\sqrt{n}\right\}  \ge \P\left\{ 
      L_n(\d_n) >0 \right\} \ge \P\left\{ L_n(\d_n)>0 ~,~ 
      n\ge \gamma\right\}. 
   \end{equation} 
   This is  at least $A_{\ref{eq:L2L}}^{-1}\d_n z_n^2 \bar\Phi(z_n) 
   \P\{n\ge \gamma\}$. 
   Therefore, the proposition follows for all $n$ large, and 
   hence all $n$ by adjusting the constants. 
\end{proof} 
 
\begin{proof}[Proof of Theorem~\ref{thm:genest}: Upper Bound] 
   Let $k=\lfloor z_n^2 \rfloor +1$, and recall the intervals 
   $I_{j,k}= [ j/k,(j+1)/k)$ for $0\le j\le k$. Then, 
   \begin{equation}\begin{split} 
      \P\left\{ \sup_{t\in E} S_n(t) \ge z_n\sqrt{n} \right\} 
         & \le \sum_{\scriptstyle 0\le j\le k:\atop 
         \scriptstyle I_{j,k}\cap E\neq\varnothing} 
         \P\left\{ \sup_{t\in I_{j,k}} S_n(t) \ge z_n\sqrt{n} \right\}\\ 
      & = \M_k(E) \P\left\{ \sup_{t\in [0,1/k]} S_n(t) 
         \ge z_n\sqrt{n}\right\}. 
   \end{split}\end{equation} 
   The last line follows from stationarity. Because 
   $\liminf_{n\to\infty} k^{-1}z_n^2=1>0$, Proposition~\ref{pr:TechnicalEst} 
   applies, and we obtain the following: 
   \begin{equation} 
      \P\left\{ \sup_{t\in E} S_n(t) \ge z_n\sqrt{n} \right\} 
      \le A_{\ref{pr:TechnicalEst}} \frac{z_n^2}{k} \M_k(E)\bar\Phi(z_n). 
   \end{equation} 
   As $n\to\infty$, $z_n^2=O(k)$, and $\M_k(E)\le 3\K_E(1/k) 
   \le 18 \K_E(z_n^{-2})$; cf.\@ Proposition~\ref{pr:KM}, 
   as well as equation~(\ref{eq:K=RV}). The probability 
   upper bound of Theorem~\ref{thm:genest} follows from this 
   discussion. 
\end{proof} 
 
\begin{proof}[Proof of Theorem~\ref{thm:genest}: Lower Bound] 
  It is likely that one can use Proposition~\ref{pr:TechnicalEst} 
  for this bound as well, but we favor a more direct approach. 
  Let $k=\K_E((16 z_n)^{-2})$, and based on this find and fix maximal 
  Kolmogorov points $x_1,\ldots,x_k$ in $E$ such that 
  whenever $i\neq j$, $|x_i-x_j|\ge (16 z_n)^{-2}$. 
  Without loss of generality, 
  we may assume that $x_1 < x_2 < \cdots < x_k$. 
  In terms of these maximal Kolmogorov points, we define 
  \begin{equation} 
     V_n = \sum_{j=1}^k \mathbf{1}_{\{ S_n(x_j) \ge 
     z_n \sqrt{n} \}}. 
  \end{equation} 
  Evidently, $\P$-almost surely, 
  \begin{equation}\label{eq:EV} 
     \EN [V_n] = k\bar\Phi(z_n) \ge \K_E\left( z_n^{-2} \right)\bar\Phi(z_n). 
  \end{equation} 
  Now we estimate the quenched second moment of $V_n$: There exists 
  an $\mathscr{N}$-measurable $\P$-almost 
  surely finite random variable $\s$ such that for all 
  $n\ge \s$, 
  \begin{equation}\begin{split} 
     \EN\left[ V_n^2 \right] &\le 2\mathop{\sum\sum}\limits_{%
        1\le i\le j\le k} \PN\left\{ S_n(x_i) \ge z_n \sqrt{n} 
        ~,~ S_n(x_j) \ge z_n \sqrt{n} \right\}\\ 
     & \le 4 \mathop{\sum\sum}\limits_{%
        1\le i\le j\le k} \exp\left( - \frac{1}{8} z_n^2 (x_j-x_i) 
        \right)\bar\Phi(z_n). 
  \end{split}\end{equation} 
  See~\ocite{KLM}*{Lemma 6.2} for the requisite joint-probability 
  estimate. 
  Whenever $j>i$, we have $x_j-x_i = \sum_{l=i}^{j-1} 
  (x_{l+1}-x_l) \ge \frac{1}{16}(j-i) z_n^{-2}$. Therefore, for all 
  $n\ge \s$, 
  \begin{equation}\label{eq:EV^2}\begin{split} 
     \EN\left[ V_n^2 \right] 
        & \le 4 \mathop{\sum\sum}\limits_{%
        1\le i\le j\le k} \exp\left( - \frac{1}{128} (j-i) 
        \right)\bar\Phi(z_n)\\ 
     & \le \frac{4}{1- e^{-1/128}} k \bar\Phi(z_n) 
        = A_{\ref{eq:EV^2}} \K_E\left( \frac{1}{16 z_n^2} \right) 
        \bar\Phi(z_n)\\ 
     & \le 6^4 A_{\ref{eq:EV^2}} \K_E\left( \frac{1}{z_n^2} \right) 
        \bar\Phi(z_n). 
  \end{split}\end{equation} 
  The last line relies on four successive applications of 
  (\ref{eq:K=RV}), and is valid if $n$ is at least 
  $r=\inf\{ k:\ z_k^2\ge 4\}$. We combine (\ref{eq:EV}), 
  (\ref{eq:EV^2}), and the Paley--Zygmund inequality to 
  deduce that for all $n\ge \s\vee r$, 
  \begin{equation} 
     \PN\left\{ V_n >0 \right\}  \ge \frac{\left( 
     \EN V_n \right)^2}{\EN\left[ V_n^2 \right]} 
     \ge \frac{1}{6^4 A_{\ref{eq:EV^2}}} 
     \K_E\left( z_n^{-2} \right)\bar\Phi(z_n), 
  \end{equation} 
  $\P$-almost surely. But for all $n\ge r$, 
  \begin{equation}\begin{split} 
     \P\left\{ \sup_{t\in E} S_n(t) \ge z_n\sqrt{n} \right\} 
        & \ge \P\left\{ V_n > 0 \right\} \ge 
        \P\left\{ V_n >0 ~,~ n\ge\s \right\}\\ 
     & \ge \frac{1}{6^4 A_{\ref{eq:EV^2}}} 
     \K_E\left( z_n^{-2} \right)\bar\Phi(z_n) \P\{ n\ge \s\}. 
  \end{split}\end{equation} 
  Because $\s$ is finite $\P$-almost surely, the lower bound 
  in Theorem~\ref{thm:genest} follows for all sufficiently large $n$, 
  and hence for all $n$ after adjusting the constants. 
\end{proof} 
 
\section{Proofs of Theorems~\ref{thm:inttest1},%
   ~\ref{thm:inttest2}, and~\ref{thm:Hdim}, 
   and Remarks%
   ~\ref{rem:inttest1} and~\ref{rem:inttest2}} 
 
The critical result is Theorem~\ref{thm:inttest1}, 
and has a long and laborious proof. Fortunately, most of this 
argument appears, in a simplified setting, in~\ocite{KLM} 
from which we borrow liberally. 
 
Throughout the following derivation, 
$\ee_n = \ee(n) = \lfloor e^{n/\ln_+(n)}\rfloor$, which is 
the so-called \emph{Erd\H{o}s sequence}. 
 
\begin{proof}[Proof of Theorem~\ref{thm:inttest1}] 
   Without loss of generality, we can assume that 
   \begin{equation}\label{eq:wlog} 
      \sqrt{\ln_+\ln_+ t} \le H(t) \le 
      2\sqrt{\ln_+\ln_+ t} 
      \qquad{}^\forall t>0. 
   \end{equation} 
   For the argument, follows~\ocite{erdos}*{eq.'s 
   (1.2) and (3.4)}. 
 
   We first dispose of the simple case $\psi_H(E)<\infty$. 
 
   By the reflection principle and by Theorem~\ref{thm:genest}, 
   \begin{equation}\begin{split} 
      \P\left\{ \max_{1\le k\le \ee(n+1)} \sup_{t\in E} S_k(t) 
         \ge H(\ee_n)\sqrt{\ee_n}\right\} & \le 2 
         \P\left\{ \sup_{t\in E} S_{\ee(n+1)}(t) 
         \ge H(\ee_n)\sqrt{\ee_n}\right\}\\ 
      & \le 2A_{\ref{eq:genest}} \K_E\left( \frac{1}{H^2(\ee_n)} \right) 
         \bar\Phi(H(\ee_n)). 
   \end{split}\end{equation} 
   Under (\ref{eq:wlog}), $\psi_H(E)$ is finite if and only if 
   $\sum_n \K_E(1/H^2(\ee_n))\bar\Phi(H(\ee_n))<\infty$. Hence, 
   the case $\psi_H(E)<\infty$ follows from a monotonicity argument. 
 
   In the case $\psi_H(E)=\infty$, define for a fixed $\vartheta>0$, 
   \begin{alignat}{2} 
      S_n^* & = \sup_{t\in E} S_{\ee(n)} (t),&\qquad 
         H_n & = H(\ee_n),\nonumber\\ 
      \mathscr{I}_n & = \left( H_n\sqrt{\ee_n}, \left( 
         H_n+\frac{\vartheta}{H_n} \right) \sqrt{\ee_n} \right],&\qquad 
         L_n & = \sum_{j=1}^n \mathbf{1}_{\{ S_j^* \in\mathscr{I}_j \}},\\ 
      f(z) & = \K_E(1/z^2)\bar\Phi(z).\nonumber 
   \end{alignat} 
   These are the present article's replacement of~\ocite{KLM}*{eq.\@ 8.10}. 
   We can choose $\vartheta$ large enough (though independent of $n$) 
   such that there exists $\eta\in(0,1)$ with the property that for 
   all $n\ge 1$, 
   \begin{equation}\label{eq:localize} 
      \eta \le \frac{ \P \left\{ S_n^*\in\mathscr{I}_n\right\}}{%
      \P\left\{ S_n^* \ge H_n\sqrt{\ee_n} \right\}} \le \eta^{-1}. 
   \end{equation} 
   To see why this holds, we mimic the proof of~\ocite{KLM}*{Lemma 8.3}, 
   but in place of their Theorem 1.4, we use Theorem~\ref{thm:genest} 
   of the present paper. 
 
   Now in light of (\ref{eq:localize}) and condition $\psi_H(E)=\infty$, 
   $\lim_{n\to\infty} \E [L_n] = \infty$. Therefore, by the Borel--Cantelli 
   lemma, it suffices to show that 
   \begin{equation}\label{eq:biggoal} 
      \limsup_{n\to\infty} \frac{ \E\left[ L_n^2 \right]}{ 
      \left( \E[ L_n ]\right)^2} <\infty. 
   \end{equation} 
   Everything comes down to estimating the following joint probability: 
   \begin{equation} 
      \mathscr{P}_{i,j} = \P \left\{ 
      S_i^* \in\mathscr{I}_i ~,~ S_j^* \in\mathscr{I}_j 
      \right\},\qquad{}^\forall j>i\ge 1. 
   \end{equation} 
   This painful task is performed by considering 
   $\mathscr{P}_{i,j}$ on three different scales: 
   (a) $j \ge i+\ln_+^{10}(i)$; 
   (b) $j \in [i+\ln_+(i) , i + \ln_+^{10}(i))$; 
   and (c) $j\in (i,i+\ln_+(i))$. 
   Fortunately, Lemmas 8.4--8.7 of~\ocite{KLM} do this for 
   us at no cost. However, we note that they hold only after we 
   replace their $S_i^*$ with ours and all multiplicative constants 
   are adjusted. Moreover, everywhere in their 
   proofs, replace ``$\sup_{t\in[0,1]}$'' by ``$\sup_{t\in E}$.'' 
   Equation~(\ref{eq:biggoal}) follows from these estimates. 
\end{proof} 
 
\begin{proof}[Proof of Theorem~\ref{thm:inttest2}] 
   First, let us suppose that $\Psi_H(E)<\infty$. 
   Then, we can write $E=\cup_{m=1}^\infty E_m$, 
   with $E_m$'s closed, such that for all $m$, $\psi_H(E_m)<\infty$. 
   Theorem~\ref{thm:inttest1} 
   proves, then, that for all $m$, 
   \begin{equation} 
      \sup_{t\in E_m} \limsup_{n\to\infty} 
      \left[ S_n(t) - H(n)\sqrt{n}\right] \le 0,\text{ a.s.} 
   \end{equation} 
   Maximize over $m=1,2,\ldots$ to prove half of Theorem~\ref{thm:inttest2}. 
 
   To prove the second half of the theorem, we assume that 
   $\Psi_H(E)=\infty$. By Lemma~\ref{lem:Psi}, 
   we can find a compact set $G\subseteq E$ such that 
   whenever $I$ is a rational interval that intersects $G$, 
   $\psi_H(I\cap G)$ is infinite. Now consider the 
   random sets 
   \begin{equation} 
      \Lambda_H^n = \left\{ t\in [0,1]:\ \sup_{\e>0}\inf_{_{\scriptstyle 
      t-\e<s<t+\e}} 
      \left[ S_n(s) - H(n)\sqrt{n} \right] >0 \right\}. 
   \end{equation} 
   By the regularity of the paths of $S_n$, 
   $\Lambda_H^n$ is open for every $n$. 
 
   By Theorem~\ref{thm:inttest2}, for any rational interval $I$ that 
   intersects $G$, 
   $\Lambda_H^n\cap (I\cap G)$ is non-empty infinitely often. In particular, 
   $\cup_{i=n}^\infty \Lambda_H^i$ intersects $I\cap G$ 
   infinitely often. 
   Therefore, we have shown that 
   $\cup_{i=n}^\infty \Lambda_H^i \cap G$ is an everywhere-dense 
   relatively-open subset of the complete compact separable metric space $G$. 
   By the Baire category theorem, 
   $\cap_{n=1}^\infty \cup_{i=n}^\infty \Lambda_H^i\cap G$ 
   is non-empty. In particular, 
   there exist uncountably-many times $t\in G\subseteq E$ 
   such that $t\in\limsup_n \Lambda_H^n = \Lambda_H$, whence the theorem. 
\end{proof} 
 
\begin{proof}[Proof of Theorem~\ref{thm:Hdim}] 
   We use a codimension argument. 
   Let $Y_\alpha$ be the stable process of \S\ref{subsec:stable} 
   which is chosen to be independent of the entire dynamical Gaussian 
   walk, and let $\mathscr{R}_\alpha$ denote its (closed) range. 
 
   By Theorem~\ref{thm:inttest1} the following are equivalent 
   for any dyadic interval $I$: 
   \begin{equation}\begin{split} 
      \limsup_{n\to\infty} 
         \left[ \sup_{t\in \mathscr{R}_\alpha\cap I} 
         S_n(t) - H(n) \sqrt{n} \right] & >0 \ 
         \Longleftrightarrow\ 
         \psi_H \left( \mathscr{R}_\alpha \cap I \right)=\infty\\ 
      \limsup_{n\to\infty} 
         \left[ \sup_{t\in \mathscr{R}_\alpha\cap I} 
         S_n(t) - H(n) \sqrt{n} \right] & \le 0 \ 
         \Longleftrightarrow\ 
         \psi_H \left( \mathscr{R}_\alpha \cap I \right)<\infty. 
   \end{split}\end{equation} 
   Recall (\ref{eq:zeta}). Thanks to (\ref{eq:stableindex}), 
   \begin{equation} 
      \E\left[ \psi_H \left(\mathscr{R}_\alpha\cap I \right)\right]\ 
      \asymp\ \int_1^\infty H^{2(1+\alpha)}(t)\bar\Phi(H(t)) \, 
      \frac{dt}{t} = \mathscr{J}_{2(1+\alpha)}(H). 
   \end{equation} 
   where `$\alpha\asymp\beta$' stands for 
   `$\alpha$ is finite if and only if $\beta$ is'. 
   Therefore, by (\ref{eq:stableindex2}) and the Paley--Zygmund 
   inequality, $\psi_H \left(\mathscr{R}_\alpha\cap I \right)$ is infinite 
   with positive probability if and only if its expectation is infinite. 
   In particular, 
   \begin{equation} 
      \P\left\{ \limsup_{n\to\infty} \sup_{t\in \mathscr{R}_\alpha\cap I} 
      \left[ S_n(t) - H(n) \sqrt{n} \right]  >0\right\}>0 \ 
      \Longleftrightarrow\ 
      \mathscr{J}_{2(1+\alpha)}(H)=\infty. 
   \end{equation} 
   Because the condition on $\mathscr{J}_{2(1+\alpha)}$ does not 
   involve the dyadic interval $I$, and since there are countably-many 
   dyadic intervals, it follows from the category portion 
   of the proof of Theorem~\ref{thm:inttest2} that 
   \begin{equation} 
      \P\left\{ \sup_{t\in \mathscr{R}_\alpha} \limsup_{n\to\infty} 
      \left[ S_n(t) - H(n) \sqrt{n} \right]  >0 \right\}>0 \ 
      \Longleftrightarrow\ 
      \mathscr{J}_{2(1+\alpha)}(H)=\infty. 
   \end{equation} 
   That is, $\Lambda_H $ intersects $\mathscr{R}_\alpha$ 
   with positive probability if and only if 
   $\mathscr{J}_{2(1+\alpha)}(H)=\infty$. 
   But it is known that $\mathscr{R}_\alpha$ can hit a set $E$ if and only if 
   $E$ has positive $(1-\alpha)$-dimensional Riesz capacity 
   $\mathrm{Cap}_{1-\alpha}(E)$~\cite{Khoshnevisan}*{Theorem 3.4.1, 
   p.\@ 384}. Thus, by the Fubini--Tonneli theorem, 
   \begin{equation} 
      \E\left[ \mathrm{Cap}_{1-\alpha}\left( \Lambda_H 
      \right) \right ]>0 \ 
      \Longleftrightarrow \mathscr{J}_{2(1+\alpha)}(H)=\infty. 
   \end{equation} 
   Because $\alpha\in(0,1)$ is arbitrary, we have shown that 
   for any $\zeta\in (0,1)$, 
   \begin{equation} 
      \E\left[ \mathrm{Cap}_{2-(\zeta/2)}\left( \Lambda_H 
      \right) \right ]>0 \ 
      \Longleftrightarrow \mathscr{J}_\zeta (H)=\infty. 
   \end{equation} 
   Frostman's theorem~\cite{Khoshnevisan}*{Theorem 2.2.1, p.\@ 521} 
   then implies the result. 
\end{proof} 
 
\begin{proof}[Proof of Remark~\ref{rem:inttest1}] 
   Because $\bar\Phi(x)\sim (2\pi)^{-1/2} x^{-1}\exp(-x^2/2)$ 
   as $x\to\infty$, 
   \begin{equation} 
      \psi_H(E)<\infty\ \Longleftrightarrow\ 
      \int_1^\infty H(t) \K_E\left( \frac{1}{H^2(t)} \right) 
      e^{-\frac12 H^2(t)}\, dt<\infty. 
   \end{equation} 
   Therefore, we can appeal to (\ref{eq:K=RV}) to see, after one or 
   two lines of calculations, that 
   \begin{equation} 
      \psi_H(E)<\infty\ \Longleftrightarrow\ 
      {}^\forall c\in\R:\ \psi_{H+(c/H)}(E)<\infty. 
   \end{equation} 
   Now we can prove the remark. 
 
   If $\psi_H(E)<\infty$, then the preceding 
   remarks and Theorem~\ref{thm:inttest1} 
   together prove that for any $c<0$, 
   \begin{equation} 
      \limsup_{n\to\infty}\left[ \sup_{t\in E} S_n(t) - \sqrt{n}\left( 
      H(n) + \frac{c}{H(n)} \right) \right] \le 0,\qquad\text{ a.s.} 
   \end{equation} 
   Thanks to (\ref{eq:wlog}), $H(n)=o(\sqrt{n})$ as $n\to\infty$. Thus, 
   let $c\to-\infty$ to see that 
   \begin{equation} 
      \limsup_{n\to\infty}\left[ \sup_{t\in E} S_n(t) - 
      H(n)\sqrt{n}\right] =-\infty,\qquad\text{ a.s.} 
   \end{equation} 
   If $\psi_H(E) = \infty$, then we argue as above, 
   but, this time, we let $c$ tend to $\infty$. 
\end{proof} 
 
\begin{proof}[Proof of Remark~\ref{rem:inttest2}] 
   We follow the proof of Remark~\ref{rem:inttest1} 
   verbatim, but apply Theorem~\ref{thm:inttest2} in place 
   of Theorem~\ref{thm:inttest1} everywhere. 
\end{proof} 

\section{Proof of Theorem~\ref{thm:invariance}} 
 
A key idea of our proof of Theorem~\ref{thm:invariance} 
is to appeal to martingale problems via 
the semi-martingale weak-convergence 
theory of~\ocite{jacod}. 
To elaborate on this connection a bit further let us note 
that $\{X_k\}_{k=1}^\infty$ are i.i.d.~copies 
of a pure-jump Feller process with generator 
\begin{equation} 
   Af(x) = \int_{-\infty}^\infty f(z)\, \nu(dz) 
   - f(x)\qquad{}^\forall f\in \mathscr{C}_0(\R). 
\end{equation} 
 
Before citing the  result of~\ocite{jacod} we need to introduce 
some more notation. This will be done in the first subsection. Let 
us note in advance that ours differs slightly from the notation 
of~\ocite{jacod}. In particular, our $B$ corresponds to their 
$B^\prime$ and our $C$ corresponds to their $\tilde{C}^\prime$. 
 
Throughout, we use the following particular construction of 
the process $\mathscr{U}$: Let $\{\beta(s,t)\}_{s,t\ge 0}$ 
denote the Brownian sheet, and define 
\begin{equation}\label{eq:UB} 
   \mathscr{U}(s,t) = \frac{\beta\left( s,e^{2t} \right)}{ 
   e^t}\qquad{}^\forall s,t\ge 0. 
\end{equation} 
The reader can check that $\mathscr{U}$ 
is indeed a continuous centered Gaussian process 
whose correlation function is given by 
(\ref{eq:covarianceU}). 
 
We aim to prove the following: 
\begin{proposition}\label{pr:limit_in_t}
   Assume, in addition, that there exists $\e>0$
   such that
   \begin{equation}\label{eq:2+e}
	   \int_{-\infty}^\infty |x|^{2+\e}\, \nu(dx)<\infty.
   \end{equation}
   Then, for each fixed $u\ge 0$, 
   $\mathscr{S}_n(u,\cdot) \Rightarrow 
   \mathscr{U}(u,\cdot)$ in the sense of 
   $\mathscr{D}([0,1])$. 
\end{proposition} 
O.\@ Rusakov~\ycite{rusakov:95}*{Theorem 3.1} 
has demonstrated that a similar
result holds for a closely-related model.
 
Because $u\mapsto \mathscr{S}_n(u,\bullet)$ is 
an infinite-dimensional L\'evy process on $\mathscr{D}([0,1])$, 
a standard argument then yields the following. 
[See Lemma 2.4 of~\ocite{eisenbaum}, 
but replace $\mathscr{D}_T(\mathscr{C}(K))$ there by 
$\mathscr{D}(\mathscr{D}([0,1]))$.] 
 
\begin{proposition}\label{prop:fdd} 
   Under the additional constraint
   \eqref{eq:2+e}, the finite-dimensional distributions of 
   $\mathscr{S}_n$ converge to those of $\mathscr{U}$. 
\end{proposition} 
In light of this, Proposition~\ref{pr:limit_in_t} and 
``tightness'' together would yield 
Theorem~\ref{thm:invariance} under \eqref{eq:2+e}. 
A truncation argument then removes \eqref{eq:2+e}.
Our proof of 
Proposition uses the machinery of~\ocite{jacod}. 
Then we follow the general outline of~\ocite{KLM}*{\S4} 
to establish tightness. 
 
\subsection{Background on Semi-Martingales} 
 
Let $\{X_t\}_{t\ge 0}$ be a cadlag semimartingale. We 
assume that $X$ is defined on the canonical 
sample space $\mathscr{D}(\R_+)$. 
 
Given a measurable function $g$, 
$\{v_t (g)\}_{t\ge 0}$ denotes the compensator of 
the process $t\mapsto 
\sum_{s \leq t,\, \Delta X_s \neq 0} g(\Delta X_s)$, 
where $\Delta X_t=X_t-X_{t-}$ designates the size of 
the jump of $X$ at time $t$. We specialize our discussion 
further by considering the subclass of processes $X$ that satisfy: 
\begin{enumerate} 
   \item $X = M + B$, where 
      $B$ is continuous and adapted, 
      and $M$ is a local martingale. 
   \item $v_t (x^2) < \infty$ for all $t$. Of course, 
      $v_t(x^2)$ stands for $v_t(g)$ where $g(x)=x^2$. 
\end{enumerate} 
For such a process $X$, write 
\begin{equation} \label{eq:C_tilde_prime_defn} 
   C_t  =  \langle M^c \rangle_t + 
   v_t (x^2) \qquad {}^\forall t\ge 0, 
\end{equation} 
where $M^c$ is the continuous part of $M$, 
and $\langle\cdot\rangle$ denotes quadratic 
variation. 
 
Let $\fclass$ denote the class of functions 
which are bounded and vanish near $0$. 
Define 
\begin{equation} \label{eq:tau_definition} 
  \tau_a  =  \inf\{ t>0:\,   |X_t| \vee |X_{t-}| \geq a \}\ 
  \qquad{}^\forall a>0. 
\end{equation} 
 
Now let $\{X^n\}_{n=1}^\infty$ denote a sequence of 
such semimartingales; 
$B^n$, $C^n$, $\tau^n_a$, 
and $v^n(g)$ denote the corresponding characteristics 
for the process $X^n$. 
 
\begin{theorem}[{\cite{jacod}*{Theorem IX.3.48}}] \label{thm:JS} 
   If the following hold for a dense subset $D$ of $\R_+$, 
   then $X_n\Rightarrow X$ in the sense of $\mathscr{D}(\R_+)$: 
   \begin{enumerate} 
      \item \label{cond:strong_major} 
         For each $a>0$ there is an increasing and 
         continuous non-random function $F^a$ 
         so that $F^a(t) - V_{\tau_a \wedge t}(B)$, 
         $F^a(t) - \langle M^c \rangle_{\tau_a 
         \wedge t}$, and $F^a(t) - v_{\tau_a \wedge t}(x^2)$ 
         are increasing functions of $t$, where $V_t(B)$ denotes 
         the total variation $B$ on $[0,t]$. 
      \item \label{cond:big_jumps} 
         For all $a>0$ and $t>0$, 
         \begin{equation} 
            \lim_{\scriptstyle 
            b \uparrow \infty} \sup_{\omega \in \mathscr{D}(\R_+)} 
            v_{\tau_a \wedge t} 
            \left(x^2 \mathbf{1}_{\{|x|>b\}} \right) 
            (\omega) =0. 
         \end{equation} 
      \item \label{cond:mg_prob_unique} 
         The martingale problem for $X$ has local uniqueness 
         in the sense of~\ocite{jacod}. 
      \item \label{cond:continuity} 
         For all $t \in D$ and $g \in \fclass$, the function 
         $\omega \mapsto (B_t(\omega),C_t(\omega), 
         v_t(g)(\omega))$ is Skorohod continuous. 
      \item \label{cond:initial} 
         $X^n_0$ converges in distribution to $X_0$. 
      \item \label{cond:fclass_vanish} 
         For all $g \in \fclass$, 
         $v^n_{t \wedge \tau^n_a} (g) - 
         v_{t \wedge \tau_a} (g) 
         \inprob 0$. 
      \item \label{cond:B_conv} 
         For all $a,t > 0$, 
         $\sup_{s \leq t} | B^n_{s \wedge \tau_a^n} 
         - B_{s \wedge \tau_a}(X^n)| \inprob 0$. 
      \item \label{cond:C_conv} 
         For all $t \in D$ and $a > 0$, 
         $C^n_{t \wedge \tau_a^n} 
         - C_{t \wedge \tau_a}(X^n) 
         \inprob 0$. 
      \item \label{cond:final} 
         For all $a,t,\e>0$, 
         \begin{equation} 
            \lim_{b \uparrow \infty} \limsup_{n \rightarrow \infty} 
            \P \left\{ {v^n_{\tau_a^n \wedge t} 
            \left( x^2\mathbf{1}_{\{|x|>b\}} \right) 
            > \e} \right\} = 0. 
         \end{equation} 
   \end{enumerate} 
\end{theorem} 
 
\subsection{Proof of Proposition~\ref{pr:limit_in_t}} 
\label{subsec:limit_in_t} 
Write $O^u_t=\mathscr{U}(u,t)$. 
We then begin by noting the 
semi-martingale characteristics of the process 
$\{O^u_t \}_{t\ge 0}$. 
First, $O^u$ solves the s.d.e., 
\begin{equation} \label{eq:OU_SDE} 
   dX_t = - X_t \,dt + \sqrt{2u} \,d\beta^u_t, 
\end{equation} 
where $\{\beta^u_t\}_{t\ge 0}$ is the 
Brownian motion $\{\beta(u,t)\}_{t \ge 0}$. 
It follows that $O^u_t = B_t(O^u) +\text{a martingale}$, where 
$B_t:\mathscr{D}(\R_+) \rightarrow \R$ is defined 
by $B_t(\omega) = -\int_0^t \omega(s)\, ds$. Also note 
that $\langle O^u \rangle_t = 2ut$.  Since 
$\{O^u_t\}_{t \geq 0}$ is path-continuous, 
$v_t(g) \equiv 0$. 
 
\begin{proof}[Proof of Proposition~\ref{pr:limit_in_t}] 
   We will verify the conditions of Theorem \ref{thm:JS} 
   as they apply to $\{S_n(u,\cdot)\}_{n=1}^\infty$ and 
   $O^u$. 
 
   The total variation of $r \mapsto B_r(\omega) = -\int_0^r \omega(s)\, ds$ 
   on $[0,t]$ is 
   $V_t(B(\omega)) = \int_0^t |\omega(s)| \, ds$. Therefore, 
   \begin{equation} 
      V_{\tau_a(\omega) \wedge t} (B(\omega)) 
      \leq a (\tau_a(\omega)\wedge t). 
   \end{equation} 
   Since $\langle M^c \rangle_t = 2ut$ and 
   $v \equiv 0$, $F^a(t)  =  [(2u \vee a)+1]t$ 
   satisfies condition \eqref{cond:strong_major}. 
 
   Condition \eqref{cond:big_jumps} is met automatically 
   because $v_t(g) \equiv 0$. 
 
   $O^{u}$ is a Feller 
   diffusion with infinitesimal drift $a(x) = -x$ 
   and infinitesimal variance $\sigma^2(x) = 2u$. In particular, $a$ 
   is Lipshitz-continuous and $\sigma^2$ is bounded. 
   Hence, by Theorems III.2.32, III.2.33, and III.2.40 
   of~\ocite{jacod}, condition \eqref{cond:mg_prob_unique} 
   is satisfied. 
 
   Because $\langle O^u \rangle_t = 2ut$, it follows that 
   $C_t = 2ut$; cf.~\eqref{eq:C_tilde_prime_defn}. 
   In particular, $\mathscr{D}(\R_+)\ni\omega\mapsto C(\omega)$ is constant. 
   Because $v_t (g)= 0$ also, 
   this establishes the continuity 
   condition \eqref{cond:continuity} for both $C$ and $v$. 
   Since 
   $\omega \mapsto \int_0^t \omega(s)\, ds$ 
   is Skorohod-continuous condition (4) is satisfied. 
 
   Condition \eqref{cond:initial} follows from 
   Donsker's Theorem; see, for example, 
   \cite{billingsley}*{Theorem 10.1}. 
 
   Fix a non-negative $g \in \fclass$, define $L = \sup_x g(x)$, and 
   suppose that $g$ vanishes on $[-\delta,\delta]$. Then, 
   we have, in differential notation, 
   \begin{equation}\begin{split} 
      d v_t^n(g) & =  \E\left[ 
         \left. d_t g \left( \Delta_t \mathscr{S}_n(u,t) 
         \right) \, \right| \, \F_t^n \right] 
         =  \sum_{k=1}^{\lfloor u n \rfloor} 
         \int_{-\infty}^\infty g \left( 
         \frac{x - X_k(t)}{\sqrt{n}} \right)\, \nu(dx)\, dt \\ 
      & \leq  L\sum_{k=1}^{\lfloor u n \rfloor} 
         \nu \{ |x - X_k(t)| \geq \sqrt{n} \delta\}\, dt \\ 
      & \leq L \sum_{k=1}^{\lfloor u n \rfloor} \frac{1}{\delta^{2+\e} 
         n^{1+(\e/2)}} 
         \int_{-\infty}^\infty |x - X_k(t)|^{2+\e}\, \nu(dx)\, dt \\ 
      & \leq \frac{2^{2+\e}  L}{\delta^{2+\e} 
         n^{1+(\e/ 2)}} \sum_{k=1}^{\lfloor u n \rfloor} 
         (\nu_{2+\e}  + |X_k(t)|^{2+\e} )\, dt. 
   \end{split}\end{equation} 
   Here, $\nu_\alpha$ denote the $\alpha$th 
   absolute moment of the measure $\nu$. 
   This implies condition~\eqref{cond:fclass_vanish}. 
 
   Next, let $\F_t=\vee_{n=1}^\infty \F^n_t$ denote the
   $\s$-algebra generated by $\{S_n(u); 0\le u\le t\}_{n=1}^\infty$,
   and note that 
   \begin{equation} 
      \E\left[  \left. 
      dX_k(t) \,\right| \, \F_t\right] 
      = \left( \int_{-\infty}^\infty x \, \nu(dx) - X_k(t) 
      \right)\, dt = - X_k(t) \, dt. 
   \end{equation} 
   Summing over $k$ gives 
   \begin{equation} 
      \E\left[ \left. d_t \mathscr{S}_n(u,t) \, 
      \right|\, \F_t \right] = - \mathscr{S}_n(u,t) \, dt. 
   \end{equation} 
   Consequently, $\mathscr{S}_n(u,t)$ has the following 
   semi-martingale decomposition: 
   \begin{equation} \label{eq:S_decomp} 
      \mathscr{S}_n(u,t) = 
      - \int_0^t \mathscr{S}_n(u,s) \,ds + \text{ a local $\F$-martingale} 
      \qquad{}^\forall t\ge 0. 
   \end{equation} 
 
   Likewise, from \eqref{eq:OU_SDE} we conclude that 
   \begin{equation} \label{eq:O_decomp} 
      O^u_t = - \int_0^t O^u_s \, ds + \text{ a local $\F$-martingale} 
      \qquad{}^\forall t\ge 0. 
   \end{equation} 
   Together \eqref{eq:S_decomp} and \eqref{eq:O_decomp} 
   verify condition \eqref{cond:B_conv}. 
 
   Because $\nu$ has mean zero and variance one, 
   \begin{equation}\begin{split} 
      d_t v_t^n (x^2) & = \E 
         \left[ \left. d_t \mathscr{S}^2_n(u,t) 
         \,\right|\, \F_t \right] \\ 
      & = \sum_{k=1}^{{\lfloor u n \rfloor}} \int_{-\infty}^\infty 
         \left( 
         \frac{x - X_k(t)}{\sqrt{n}} \right)^2\, 
         \nu(dx) \, dt \\ 
      & = \left( \frac{{\lfloor u n \rfloor}}{n} + 
         \frac{1}{n}\sum_{k=1}^{{\lfloor u n \rfloor}} 
         (X_k(t))^2 \right) \, dt. 
   \end{split}\end{equation} 
   The pure-jump character of $\mathscr{S}_n(u,\cdot)$ 
   implies that the quadratic variation of 
   the continuous part of the local martingale in 
   \eqref{eq:S_decomp} is zero, whence 
   $C^n_t = v_t^n (x^2)$.  By the computation 
   above and the law of large numbers, 
   $C^n_t \inprob 2u t =C_t = \langle O^u \rangle_t$. 
   Therefore, condition \eqref{cond:C_conv} 
   is satisfied. 
 
   Finally, after recalling that $\nu_\alpha$ is the $\alpha^{\text{%
   \underline{th}}}$ 
   absolute moment of $\nu$, we have 
   \begin{equation}\begin{split} 
      &v_t ^n\left( x^2 \mathbf{1}_{\{ |x| > b \}} \right) 
         = \int_0^t \E\left[ \left. d_s \mathscr{S}_n^2 (u,s) 
         \mathbf{1}_{\{ |d_s \mathscr{S}_n(u,s)| > b \}} \, 
         \right| \, \F_s \right] \\ 
      & = \frac{1}{n}  \int_0^t 
         \sum_{k=1}^{\lfloor u n \rfloor} \int_{-\infty}^\infty 
         \left( x - X_k(s)\right)^2 \mathbf{1}_{\{ |x - X_k(s)| > b\}} 
         \, \nu(dx)\, ds\\ 
      & \leq \frac{1}{n}  \int_0^t 
         \sum_{k=1}^{\lfloor u n \rfloor} 
         \left[ \int_{-\infty}^\infty ( x - X_k(s))^{2+\e} 
         \, \nu(dx)\right]^{2/(2+\e) } \\ 
      & \qquad \times \left[ \nu\{ 
         |x - X_k(s)| > b\} \right]^{\e/(2+\e) }\, ds\\ 
      & \leq \frac{2^{2+\e} }{n b^{\e/(2+\e) }} \int_0^t 
         \sum_{k=1}^{\lfloor u n \rfloor} 
         \left( \nu_{2+\e}  + |X_k(s)|^{2+\e} 
         \right)^{2/(2+\e) } 
         \left( \nu_1 + |X_k(s)|\right)^{ \e/(2+\e) }\, ds. 
   \end{split}\end{equation} 
   By the stationarity of $X$, 
   \begin{equation}\begin{split} 
       &\E\left[ v_t^n\left( x^2 \mathbf{1}_{ 
          \{ |x| > b \}} \right)\right]\\ 
       &\leq \frac{2^{2+\e} t}{b^{\e/(2+\e) }} 
          \E \left[ \left( \nu_{2+\e}  + |X_1(0)|^{2+\e}  \right)^{2/(2+\e) } 
          \left( \nu_1 + |X_1(0)|\right)^{\e/(2+\e) } \right]. 
   \end{split}\end{equation} 
   Also, since $t \mapsto v_t^n( x^2 \mathbf{1}_{\{|x|>b\}})$ 
   is non-decreasing we have 
   \begin{equation} 
      \E \left[ v_{\tau^n_a \wedge t}^n\left( x^2 \mathbf{1}_{\{ |x|>b \}} 
      \right) \right] 
      \leq \frac{A t}{b^{\e/(2+\e)}}. 
   \end{equation} 
   Therefore, by Markov's inequality, 
   condition \eqref{cond:final} holds. 
\end{proof} 
 
\subsection{Tightness} 
This portion contains a variation on the argument 
in \ocite{KLM}*{\S4}. We appeal to a criterion 
for tightness in $\mathscr{D}([0,1]^2)$ 
due to Bickel and Wichura~\ycite{bickel}. 
[Because $\mathscr{D}([0,1]^2) \simeq
\mathscr{D}(\mathscr{D}([0,1]))$,
we will not make a distinction between the two spaces.]
 
A \emph{block} is a two-dimensional half-open rectangle 
whose sides are parallel to the axes; i.e., $I$ is a block 
if and only if it has the form $(s,t]\times(u,v]\subseteq( 0,1]^2$. 
Two blocks $I$ and $I'$ are \emph{neighboring} if either: 
(i) $I=(s,t]\times(u,v]$ and $I'=(s',t']\times(u,v]$ (horizontal 
neighboring); or 
(ii) $I=(s,t]\times(u,v]$ and $I'=(s,t]\times(u',v']$ 
(vertical neighboring). 
 
Given any two-parameter stochastic process $Y = \{ Y(s,t);\ 
s,t\in[0,1]\}$, and any block $I = (s,t]\times(u,v]$, 
the \emph{increment of $Y$ over $I$} [written as 
$\Delta Y(I)$] is defined as 
\begin{equation} 
   \Delta Y(I)  =  Y(t,v)-Y(t,u)-Y(s,v)+Y(s,u). 
\end{equation}

\begin{lemma}[{Refinement to~\ocite{bickel}*{Theorem 3}}] 
  \label{lem:BW} 
   Let $\{ Y_n\}_{n=1}^\infty$ denote a sequence of random fields 
   in $\mathscr{D}([0,1]^2)$ such 
   that for all $n\ge 1$, $Y_n(s,t)=0$ if $st=0$. Suppose 
   that there exist constants $A_{\ref{eq:BW}}>1$, 
   $\theta_1,\theta_2,\gamma_1,\gamma_2>0$ such that 
   they are all independent of $n$, and 
   whenever $I = (s,t]\times(u,v]$ and 
   $J = (s',t']\times(u',v']$ are neighboring blocks, and if 
   $s,t,s',t'\in n^{-1}\mathbf{Z} \cap[0,1]$, then 
   \begin{equation}\label{eq:BW} 
      \E\left[ \,\left| \Delta Y_n (I) 
      \right|^{\theta_1} 
      \left| \Delta Y_n(J)  \right|^{\theta_2} \,\right] \le 
      A_{\ref{eq:BW}} \left| I\right|^{\gamma_1} 
      \left| J\right|^{\gamma_2}, 
   \end{equation} 
   where $|I|$ and $|J|$ denote respectively the planar Lebesgue measures 
   of $I$ and $J$. 
   If, in addition, $\gamma_1+\gamma_2>1$, then 
   $\{ Y_n\}_{n=1}^\infty$ is a tight sequence. 
\end{lemma}

Additionally, we need the following \emph{a priori} estimate.
\begin{lemma}\label{lem:tightness}
   In Theorem~\ref{thm:invariance},
   \begin{equation}
	   \E \left[ \max_{_{\scriptstyle k\in\{1,\ldots,n\}}}
	   \sup_{u\in[0,1]} | S_k(u) |^2 \right] \le 64n
	   \qquad {}^\forall n\ge 1.
   \end{equation}
\end{lemma}

\begin{proof}
   We choose and fix an integer $n\ge 1$.
   Also, we write $\EN$ for the conditional-expectation
   operator $\E[\cdots\,|\,\mathscr{N}]$, where $\mathscr{N}$
   denotes the $\s$-algebra generated by the clocks.

   We can collect the jump-times of the
   process $\{S_i(u)\}_{u\in [0,1]}$ for all $i=1,\ldots,n$.
   These times occur at the jump-times of a homogeneous,
   mean-$n$ Poisson process by time one. Define $T_0=0$ and enumerate
   the said jumps to obtain $0=T_0<T_1<T_2<\ldots<T_{N(n)}$.
   The variable $N(n)$ has the Poisson distribution with mean $n$.
   
   If $u\in[T_j,T_{j+1})$, then $S_n(u)=S_n(T_j)
   =\sum_{\ell=0}^{j-1}\{ S_n(T_{\ell+1}) - S_n(T_\ell)\}$. 
   This proves that
   \begin{equation}
      \sup_{u\in[0,1]} \left| S_n(u) \right| 
      = \max_{1\le j\le N(n)} \left| \sum_{\ell=0}^{j-1}
      \zeta_\ell \right|.
   \end{equation}
   Here, the $\zeta$'s are independent of $\mathscr{N}$,
   and have the same distribution as $\nu\star \nu^-$ where
   $\nu^-(G)=\nu(-G)$.  Moreover, the $\zeta_{2i}$'s [resp.\@
   $\zeta_{2i+1}$'s] form an independent collection. 
   In accord with Doob's maximal $(2,2)$-inequality,
   \begin{equation}\begin{split}
      &\EN\left[ \sup_{u\in[0,1]} \left| S_n(u) \right|^2 \right]\\ 
      &\le 2\left\{ \EN\left[
         \max_{1\le j\le N(n)} \left| \sum_{\ell<j:~
         \mathrm{odd}} \zeta_\ell \right|^2 \right] 
         + \EN\left[
         \max_{1\le j\le N(n)} \left| \sum_{\ell<j:~
         \mathrm{even}} \zeta_\ell \right|^2
         \right]\right\}\\
      & \le 8 \EN\left[ \sum_{\ell=0}^{N(n)-1}
         \zeta_\ell^2 \right] = 16 N(n)\qquad
         \mathrm{a.s.}
   \end{split}\end{equation}
   [We have used also the inequality $(x+y)^2\le 2\{x^2+y^2\}$.]
   Take expectations to obtain
   \begin{equation}\label{predoob}
	   \E\left[ \sup_{u\in[0,1]} |S_n (u)|^2 \right] \le 16n.
   \end{equation}
   It is easy to see that $n\to \sup_{u\in[0,1]}|S_n(u)|$ is
   a submartingale. Thus, Doob's strong $(2,2)$-inequality
   and \eqref{predoob} together imply the lemma.
\end{proof}
 
\subsection{Proof of Theorem~\ref{thm:invariance}}
We proceed in two steps. 

\emph{Step 1. The $L^4(\P)$ Case.}\
First we derive the theorem when 
$\E\{|X_0(u)|^4\}$ is finite. In this case, \eqref{eq:2+e}
holds and so it remains to derive tightness. We do so by
appealing to Lemma~\ref{lem:BW}.

Consider first the vertical neighboring case.
By the stationarity 
of the increments of random walks
we need only consider the case where $I = (0,s] \times (0,u]$ and 
$J = (0,s] \times (u,v]$, where $s \in n^{-1} \mathbf{Z}$. Clearly,
\begin{equation}\begin{split} 
    \Delta \mathscr{S}_n(I)  &= \mathscr{S}_n(s,u) - 
       \mathscr{S}_n(s,0) = 
       \frac{S_{\floor{sn}}(u) - 
       S_{\floor{sn}}(0)}{\sqrt n},\\ 
    \Delta \mathscr{S}_n(J)  &= 
       \mathscr{S}_n(s,v) - \mathscr{S}_n(s,u) = \frac{ 
       S_{\floor{sn}}(v) - S_{\floor{sn}}(u)} 
       {\sqrt n}. 
\end{split}\end{equation} 
By the Cauchy-Schwarz inequality, 
$\|\Delta \mathscr{S}_n(I) \Delta \mathscr{S}_n(J)
\|_2^2 \le \|\Delta \mathscr{S}_n(I) 
\|_4^4 \ \|\Delta \mathscr{S}_n(J)\|_4^4$.
Note that the distribution of $\Delta \mathscr{S}_n(I)$ [resp.\@
$\Delta \mathscr{S}_n(J)$] is the same as
$\sum_{i=1}^{N_n(u)} (\xi_i-\xi_i')$
[resp.\@ $\sum_{i=1}^{N_n(v-u)} (\xi_i-\xi_i')$], where: 
(i) $\{\xi_i\}_{i=1}^\infty$
is an i.i.d.\@ sequence, each distributed according to $\nu$;
(ii) $\{\xi_i'\}_{i=1}^\infty$ is an independent copy of
$\{\xi_i\}_{i=1}^\infty$; and
(iii) $N_n(r)$ is a Poisson random variable, with mean 
$\lfloor nr\rfloor$, that is
independent of all of the $\xi$'s.
These remarks, together with a direct computation, show that
there exists a finite constant $K$ such that
$\|\Delta \mathscr{S}_n(I) \Delta \mathscr{S}_n(J)\|_2
\le K |I| \ |J|$. A similar inequality is valid for 
the horizontal neighboring case. That is simpler to derive
than the preceding, and so we omit the details.
This and Lemma~\ref{lem:BW} together prove tightness in
the case that the $X_k(0)$'s are in $L^4(\P)$. According
to Proposition~\ref{prop:fdd}, Theorem~\ref{thm:invariance}
follows suit in the case that $X_1(0)\in L^4(\P)$.
   
\emph{Step 2. Truncation.}\
Now we prove Theorem~\ref{thm:invariance}
under the conditions give there; that is,
$\int_{-\infty}^\infty x\, \nu(dx)=0$
and $\int_{-\infty}^\infty x^2\, \nu(dx)=1$.

For any $c>0$ define $X_k^c(u)=
X_k(u) \mathbf{1}_{\{ |X_k(u)|\le c\}} -
\int_{-c}^c x\, \nu(dx)$. Also define
$S_n^c(u) = \sum_{k=1}^n X_k^c(u)$. It is easy to see
that $\{S_n^c\}_{n=1}^\infty$ and $\{S_n-S_n^c\}_{n=1}^\infty$ 
define two independent, centered,
dynamical random walks.  According to Step 1,
$\s(c)\mathscr{S}_n^c\Rightarrow \mathscr{U}$ as $n\to\infty$,
where: (a) $\mathscr{S}^c$ is defined as $\mathscr{S}$,
but in terms of the $X^c$'s instead of the $X$'s; and (b)
$\s^{-2}(c) = \mathrm{Var}(X_1(0); |X_1(0)|\le c)$.
Because $\lim_{c\to\infty} \s(c)=1$ and $\mathscr{U}$
is continuous it suffices to prove that for all $\e>0$,
\begin{equation}\label{goal:tightness}
     \lim_{_{\scriptstyle c\to\infty}}
     \sup_{n\ge 1 }  \P\left\{ 
     \sup_{s,t\in[0,1]} \left|
     \mathscr{S}_n(s,t) - 
     \mathscr{S}^c_n (s,t) \right| \ge \e  \right\} =0.
\end{equation} 
But we can change scale and apply
Lemma~\ref{lem:tightness} to deduce that
\begin{equation}\
     \E\left[ 
     \max_{_{\scriptstyle k\in\{ 1,\ldots,n\}}}
     \sup_{u\in[0,1]} \left| S_k(u) -
     S^c_k(u) \right|^2 \right] \le 
     64\mathrm{Var}\left( X_1(0); |X_1(0)|\ge c \right)n,
\end{equation}
for all integers $n\ge 1$. Equation
\eqref{goal:tightness} follows from the preceding
and the Chebyshev inequality;
Theorem~\ref{thm:invariance} follows.
\hfill$\square$ 
 
\section{Proof of Theorem~\ref{thm:gamblersruin}} 
 
First, we develop some estimates for general random walks. 
Thus, for the time being, let $\{s_n\}_{n=1}^\infty$ 
denote a random walk on $\mathbf{Z}$ with increments 
$\{\xi_n\}_{n=1}^\infty$. 
As is customary, 
let $P_x$ denote the law of $\{x+s_n\}_{n=1}^\infty$ 
for any $x\in\R$, and introduce $s_0$ 
so that $P_z\{s_0=z\}=1$ for all $z\in\mathbf{Z}$; 
note that $\P=P_0$. We assume, for the time being, 
that the set of possible points of $\{s_n\}_{n=1}^\infty$ 
generates the entire additive group $\mathbf{Z}$. 
Thanks to the free abelian-group theorem this is 
a harmless assumption. 
See~\ocite{Khoshnevisan}*{p.~78} for details. Define 
\begin{equation} 
   G(n)  = \sum_{i=1}^n P_0\{ s_i =0\}\qquad 
   {}^\forall n\ge 1. 
\end{equation} 
 
\begin{lemma}\label{lem:PGP} 
   For all $n\ge 1$ and $z\in\mathbf{Z}$, 
   \begin{equation} 
      P_z\{ T(0)>n\} \le \frac{1}{G(n) P_0\{T(z) \le T(0)\}}. 
   \end{equation} 
\end{lemma} 
 
\begin{proof} 
   We start with a last-exit decomposition. Because 
   the following are disjoint events, 
   \begin{equation}\begin{split} 
      1 & \ge \sum_{j=1}^n P_0\left\{ s_j =0, 
         s_{j+1}\neq 0, \ldots, s_{j+n} \neq 0\right\} \\ 
      & = \sum_{j=1}^n P_0\left\{ s_j=0 , s_{j+1}-s_j\neq 0, 
         \ldots, s_{j+n}-s_j\neq 0 \right\}\\ 
      & = \sum_{j=1}^n P_0\{ s_j=0\} P_0\{T(0)>n\}\\
      & = G(n) P_0\{T(0)>n\}. 
   \end{split}\end{equation} 
   By the strong Markov property, 
   \begin{equation}\label{eq7} 
      P_0\{ T(0)>n\} \ge P_0\{ T(z) \le T(0)\} P_z\{T(0)>n\}. 
   \end{equation} 
   The result follows from this 
   and the preceding display. 
\end{proof} 
 
Consider the local times, 
\begin{equation} 
   L^x_n = \sum_{j=0}^n \mathbf{1}_{\{s_j=x\}}\qquad 
   {}^\forall x\in\mathbf{Z}, \, n\ge 0. 
\end{equation} 
Evidently, $G(n)= E_0[L^0_n]-1$, where 
$E_z$ denotes the expectation operator under $P_z$. 
 
\begin{lemma}\label{lem:PL} 
   For all $z\in\mathbf{Z}$ and $n\ge 1$, 
   $P_z\{T(0)>n\} \le E_0[L^0_{T(z)}]/G(n)$. 
\end{lemma} 
 
\begin{proof} 
   If $z=0$, then $L^0_{T(z)}=2$, and the lemma follows from
   Lemma~\ref{lem:PGP}. From now on, we assume that $z\neq 0$.
   We can apply the strong Markov 
   property to the return times to $z$, 
   and  deduce that for all 
   non-negative integers $k$, 
   \begin{equation}\label{eq:LT-geom} 
      P_0\left\{ L_{T(z)}^0 = k+1 \right\} 
      =\left[ P_0\left\{ T(0) < T(z) \right\} \right]^k 
      P_0\{ T(z) < T(0)\}. 
   \end{equation} 
   [The $k+1$ is accounted for by the fact that $L^0_0=1$.]
   Therefore, the $P_0$-law of $L^0_{T(z)}$ is geometric 
   with mean 
   \begin{equation}\label{eq:mean-geom} 
      E_0 \left[ L^0_{T(z)} \right] = \frac{1}{P_0 
      \left\{ T(z) < T(0) \right\}}. 
   \end{equation} 
   This and Lemma~\ref{lem:PGP} together prove the lemma. 
\end{proof} 
 
\begin{lemma}\label{lem:rec-est} 
   If $\{s_n\}_{n=1}^\infty$ is recurrent, then for all 
   non-zero integers 
   $z$ and all $n\ge 1$, 
   \begin{equation}\begin{split} 
      &P_z \left\{ T(0)>n \right\} \le \frac{ 
         2\{ 1+G(\theta(z))\}}{G(n)}, \text{ where}\\ 
      &\theta(z) = \inf\left\{ n\ge 1:\ P_0 
         \left\{ T(z) > n \right\} \le \frac18 \right\}. 
   \end{split}\end{equation} 
\end{lemma} 
 
\begin{proof} 
   Recurrence insures that $\theta(z)$ is finite for 
   all $z\in\mathbf{Z}$. Now for any positive integer $m$, 
   \begin{equation}\begin{split} 
      E_0\left[ L^0_{T(z)} \right] & \le E_0 
         \left[ L^0_m \right] + E_0\left[ 
         L^0_{T(z)} ; T(z) >m \right]\\ 
      & \le 1+G(m) + \sqrt{E_0\left[ \left( L^0_{T(z)} 
         \right)^2\right] P_0\left[ T(z)>m\right\}}. 
   \end{split}\end{equation} 
   Since $L^0_{T(z)}$ 
   has a geometric distribution [see (\ref{eq:LT-geom})], 
   $E_0[(L^0_{T(z)})^2]\le 2 \{ E_0[ L^0_{T(z)}] \}^2$. 
   Thus, 
   \begin{equation}\label{eq:G(m)} 
      E_0\left[ L^0_{T(z)} \right] \le 1+
      G(m) + E_0\left[ L^0_{T(z)}\right] 
      \sqrt{2 P_0\left[ T(z)>m\right\}}. 
   \end{equation} 
   Choose $m=\theta(z)$ to find that the square root is 
   at most $\frac12$. Solve for $E_0[ L^0_{T(z)} ]$ 
   to finish. 
\end{proof} 
 
\begin{lemma}\label{lem:hit-est} 
   Suppose $\E[\xi_1]=0$ and $\s^2=\E[\xi_1^2]<\infty$. 
   Then we can find a finite constant 
   $A_{\ref{eq:hit-est}}>1$ such that 
   \begin{equation}\label{eq:hit-est} 
      P_z\{ T(0)>n \} \le A_{\ref{eq:hit-est}} 
      \frac{1+|z|}{\sqrt n}\qquad 
      {}^\forall z\in\mathbf{Z},\, n\ge 1. 
   \end{equation} 
\end{lemma} 
 
\begin{proof} 
   First of all, we claim that there exists $A_{\ref{eq:claim1}}>1$ 
   such that for all $n\ge 1$, 
   \begin{equation}\label{eq:claim1} 
      A_{\ref{eq:claim1}}^{-1} 
      \sqrt n \le G(n) \le A_{\ref{eq:claim1}} \sqrt n. 
   \end{equation} 
   When $\{s_n\}_{n=1}^\infty$ is strongly aperiodic 
   this follows from the local
   central limit theorem~\cite{spitzer}*{II.7.P9}. In the general 
   case, consider the random walk $\{s_n'\}_{n=1}^\infty$ whose 
   increment-distribution is $\frac12(\nu+\d_0)$.
   The walk $\{s_n'\}_{n=1}^\infty$ 
   has the same law as $\{s_{c(n)}\}_{n=1}^\infty$ 
   where $c(n) =\min\{m:\ \lambda_0+\cdots+\lambda_m\ge n\}$ 
   for an i.i.d.~sequence $\{\lambda_n\}_{n=1}^\infty$ of 
   mean-$(\frac12)$ geometric random variables 
   that are totally independent of 
   $\{s_n\}_{n=1}^\infty$. Because $\sum_{i=0}^n 
   \mathbf{1}_{\{ s_i'=0\}}=\sum_{i=0}^n \lambda_i\mathbf{1}_{\{ s_i=0\}}$, 
   it follows that $G'(n)=2G(n)$ where 
   $G'(n)=\sum_{i=1}^n P\{ s_i' =0\}$. Because $\{s_n'\}_{n=1}^\infty$ is 
   strongly aperiodic, (\ref{eq:claim1}) follows. In light of this 
   and Lemmas~\ref{lem:PGP}
   and~\ref{lem:rec-est}, it suffices to prove that 
   \begin{equation}\label{eq:theta-O} 
      \theta(z) = O(z^2)\quad \text{as $|z|\to\infty$ in $\mathbf{Z}$}.
   \end{equation} 
   If $\beta>0$ is fixed, then 
   \begin{equation}\begin{split} 
      P_z \left\{ T(0)>\lfloor \beta z^2\rfloor\right\} 
         & = P_0\left\{ L^{-z}_{\lfloor \beta z^2\rfloor} =0 \right\} 
         \le P_0\left\{ L^{-z}_{\lfloor \beta z^2\rfloor} \le 
         \sqrt{|z|}\right\}\\ 
      & = P_0\left\{ \ell_{\s\beta}^{-1} \le \s \right\} + o(1) 
         \quad\text{as $|z|\to\infty$}. 
   \end{split}\end{equation} 
   Here $\ell^{-1}_t$ denotes the local time of 
   Brownian motion at $-1$ by time $t$. [The preceding display 
   follows from the local-time invariance principle of 
   \ocite{borodin}.] Recurrence of Brownian motion implies that 
   there exist $\beta,z_0>0$ such that whenever 
   $|z|\ge z_0$, $P_z\{ T(0)>\beta z^2\}\le\frac18$; 
   i.e., $\theta(z)\le\beta z^2$ as long as $|z|\ge z_0$. 
   This verifies (\ref{eq:theta-O}) and completes our proof. 
\end{proof} 
 
\begin{proof}[Proof of Theorem~\ref{thm:gamblersruin}] 
   We can appeal to the free abelian-group theorem again 
   to assume without loss of generality that the possible 
   values of $\{s_n\}_{n=1}^\infty$ generate the entire 
   additive group $\mathbf{Z}$. 
 
   Apply (\ref{eq:G(m)}) with $m=\theta(z)$ to find that 
   $E_0[L^0_{T(z)}]\le 2\{1+G(\theta(z))\}$. 
   Combine this with (\ref{eq:claim1}) and (\ref{eq:theta-O}) 
   to find that $E_0[L^0_{T(z)}] \le A\sqrt{1+z^2}$ for 
   some constant $A$ that does not depend on $z\in\mathbf{Z}$. 
   This and (\ref{eq:mean-geom}) together imply the lower bound of 
   Theorem~\ref{thm:gamblersruin}. 
 
   To obtain the other bound let $\tau=\inf\{n: s_n\le 0\}$. 
   Because $T(0)\ge \tau$, Lemma~\ref{lem:hit-est} and (\ref{eq7}) 
   together prove that 
   \begin{equation} 
      \frac{A_{\ref{eq:hit-est}}}{\sqrt n} 
      \ge P_0\{ T(z)\le T(0)\} P_z\{ \tau >n\}. 
   \end{equation} 
   Thanks to~\ocite{pemantleperes}*{Lemma 3.3}, as long as 
   $|z|\le A'\sqrt n$ for a fixed $A'$, 
   $P_z\{\tau>n\}\ge A''|z|/\sqrt n$ for a fixed 
   $A''$. The result follows. 
\end{proof} 
 
\begin{remark} 
   The last portion of the preceding proof shows also that 
   $P_z\{ T(0)>n\} \ge A''|z|/\sqrt n$. This proves that 
   the bound in Lemma~\ref{lem:hit-est} is sharp up to a 
   multiplicative constant. 
\end{remark} 
 
\section{Proof of Theorem~\ref{thm:DS-REC}} 
The basic outline of our proof 
follows the same general line of thought 
as the derivation of (3.1) of~\ocite{penrose}. 
However, as was noted by 
\ocite{benjamini}, the present discrete set-up 
contains inherent technical difficulties that do 
not arise in the continuous setting of \ocite{penrose}. 

Choose and fix a large positive integer $M$, and define 
\begin{equation}\label{eq:abg}
   \gamma = \frac{3}{6+2\e},\qquad
   q_n = \left\lfloor\frac{n}{M} \right\rfloor, 
   \qquad{}^\forall n=1,2,\ldots\,. 
\end{equation} 
Within $[n/2,n]$ we can find 
$\lfloor n/(4q_n)\rfloor $-many closed 
intervals $\{I_k^n\}_{k=1}^{\lfloor n/(4q_n)\rfloor}$, 
of length $q_n$ each, such that the distance between 
$I^n_i$ and $I^n_j$ is at least $q_n$ if $i\neq j$. 
Motivated by \S5 of 
\ocite{benjamini}, let $E_n(t)$ denote the event that 
\begin{equation}\begin{split} 
   &\{S_k(t)\}_{k=0}^\infty\text{ takes both (strictly) positive 
      and  (strictly) negative}\\ 
   &\text{values in every one of $I_1^n,\ldots, I_{\lfloor
   n/(4q_n)\rfloor}^n$}. 
\end{split}\end{equation} 
Also let $F_n(t)$ denote the event that 
\begin{equation} 
   \{S_k(t)\}_{k=0}^\infty\text{ does not return to zero 
   in $[n/2,n]$}. 
\end{equation} 
 
\begin{lemma}\label{lem:estest} 
   Uniformly for all $t\ge 0$, 
   \begin{equation} 
      \limsup_{n\to\infty} \frac{\ln \P(E_n(t)\cap F_n(t))}{\ln n}
      \le - \frac{M \gamma}{12}.
   \end{equation} 
\end{lemma} 
 
\begin{proof} 
   The uniformity assertion holds tautologically since 
   $\P(E_n(t) \cap F_n(t))$ 
   does not depend on $t\ge 0$. Without loss 
   of generality, we may and will work 
   with $t=0$. 
 
   Let $f^n_i$ denote the smallest value in $I^n_i$. 
   Also define 
   \begin{equation} 
      c^n_i = \inf \left\{ 
      \ell\in I^n_i\setminus\{f^n_i\}:\ 
      S_{\ell-1}(0)S_\ell(0)<0 \right\}, 
   \end{equation} 
   where $\inf\varnothing=\infty$. 
   Finally, define $A^n_i$ to be the event that 
   $c^n_i$ is finite, but $S_k(0)\neq 0$ 
   for all $k=c^n_i+1,\ldots,c^n_i+q_n$. A little thought 
   shows that for any integer $j\ge 1$, 
   \begin{equation}\label{AAAA}\begin{split} 
      &\P\left( \left. A^n_{j+1}\, \right| 
         \, A^n_1,\ldots,A^n_j\right)\\ 
      & \le \P \left\{ \max_{1\le i\le n} |X_i(0)|\ge n^\gamma \right\} 
         + \sup_{|x|\le n^\gamma} P_x \left\{ S_k(0) 
         \neq 0,\ {}^\forall k=1,\ldots, q_n\right\}. 
   \end{split}\end{equation}
   To estimate the first term we note that
   $(2+\e)\gamma-1=\gamma/3$. Therefore,
   \begin{equation}\begin{split}
      \P \left\{ \max_{1\le i\le n} |X_i(0)|\ge n^\gamma \right\} &
         \le n \P\left\{ |X_1(0)| \ge n^\gamma\right\}
         \le \frac{\E \left\{ |X_1(0)|^{2+\e}\right\} }{
         n^{-(2+\e)\gamma-1}}\\
      & = O\left( n^{-\gamma/3} \right)\qquad\text{ as $n\to\infty$}.
   \end{split}\end{equation}
   See~\eqref{eq:abg}. On the other hand,
   by Lemma~\ref{lem:hit-est} and~\eqref{eq:abg},
   \begin{equation}
      \sup_{|x|\le n^\gamma} P_x \left\{ S_k(0) 
      \neq 0,\ {}^\forall k=1,\ldots, q_n\right\}
      \le A_{\ref{eq:hit-est}} \frac{n^\gamma}{\sqrt{q_n}}
      = O\left( n^{-\gamma/3} \right),
   \end{equation}
   as $n\to\infty$.
   These remarks, together with~\eqref{AAAA} imply that 
   \begin{equation} 
         \sup_{j\ge 1} \P\left( \left. A^n_{j+1}\, \right| 
         \, A^n_1,\ldots,A^n_j\right)= 
         O\left( n^{-\gamma/3} \right).
   \end{equation}   
   Thus, as $n\to\infty$, 
   \begin{equation}\begin{split} 
      \P(E_n(t)\cap F_n(t)) & 
         \le \P\left( \bigcap_{i=1}^{\lfloor n/(4q_n)\rfloor} 
         A^n_i \right)\\
      & = O\left( n^{-\gamma\lfloor n/(4q_n)\rfloor
         /3}\right)  \le n^{o(1)- M\gamma/12}.
   \end{split}\end{equation}
   This proves the lemma. 
\end{proof} 
 
\begin{lemma}\label{lem:5.3_benjamini}
   There exists $M_0=M_0(\e)$ such that whenever
   $M>M_0$, 
   \begin{equation}
      \sum_{n=1}^\infty \P \left( \bigcap_{s\in[0,1]} \left[
      E_n(s)\cap F_n(s) \right] \right) <\infty.
   \end{equation} 
\end{lemma} 
 
\begin{proof} 
   By Lemma~\ref{lem:estest} 
   and the strong Markov property, 
   \begin{equation}\label{eq:5.3_benjamini1} 
      \int_0^\infty \P\left( \bigcap_{s\in[0,t]} \left[ 
      E_n(s)\cap F_n(s) \right] 
      \right) e^{-t}\, dt\le n^{o(1)+1-M\gamma/12}\qquad(n\to\infty). 
   \end{equation} 
   See the proof of Lemma 5.3 of~\ocite{benjamini}. 
   On the other hand, 
   \begin{equation}\begin{split}
      &\int_0^\infty \P\left( \bigcap_{s\in[0,t]} \left[ 
	      E_n(s)\cap F_n(s) \right] 
	      \right)e^{-t}\, dt\\
      &\ge \frac1e \int_0^1\P\left( \bigcap_{s\in[0,t]} \left[ 
	      E_n(s)\cap F_n(s) \right] 
	      \right)\,dt\\
	&\ge \frac 1e \P\left( \bigcap_{s\in[0,1]}
		\left[ E_n(s)\cap F_n(s)\right]
		\right).
   \end{split}\end{equation}
   Therefore, 
   $\P(\cap_{s\in[0,1]} [E_n(s)\cap F_n(s)])\le n^{o(1)+1-M\gamma/12}$. 
   The lemma follows with $M_0=24/\gamma$. 
\end{proof} 
 
The following is essentially Lemma 5.4 of~\ocite{benjamini}. 
To prove it, go through their derivation, and 
replace their $I^n_i$'s by ours. 
 
\begin{lemma}\label{lem:5.4_benjamini}
   Suppose $M>M_0$. Then,
   \begin{equation}
      \P \left( \bigcap_{t\ge 0} \limsup_n E_n(t) \right)=1.
   \end{equation}
\end{lemma} 
 
\begin{proof}[Proof of Theorem~\ref{thm:DS-REC}]
   Choose and fix some $M>M_0$. Then
   follow along the proof of~\ocite{benjamini}*{Theorem 1.11}, 
   but replace their $\tau$ by one, and 
   the respective applications of their 
   Lemmas 5.3 and 5.4 by our Lemmas~\ref{lem:5.3_benjamini} 
   and~\ref{lem:5.4_benjamini}. 
\end{proof} 
 
\section{Applications to the OU Process on Classical Wiener Space} 
\label{sec:OU} 
 
Let $\beta$ denote a two-parameter Brownian sheet and consider 
once more the construction \eqref{eq:UB}. In addition, 
recall from \S\ref{subsec:limit_in_t} the process 
$\{O^u_t; u\in[0,1]\}_{t\ge 0}$, which can be written in terms of 
the Brownian sheet $\beta$ as follows: 
\begin{equation} 
   O^u_t = \frac{\beta(u,e^{2t})}{e^t} 
   \qquad{}^\forall t\ge 0, u\in[0,1]. 
\end{equation} 
This proves readily that the process 
$\{O^\bullet_t\}_{t\ge 0}$ is an infinite-dimensional 
stationary diffusion on $C([0,1])$ whose invariant measure 
is the Wiener measure on $C([0,1])$. The process 
$O=\{O^\bullet_t\}_{t\ge 0}$ is a fundamental object in 
infinite-dimensional analysis. 
See, for example,~\ocite{kuelbs},~\ocite{malliavin}, 
and~\ocite{walsh}. These furnish three diverse theories 
in each of which $O$ plays a central role. 
 
An interesting artifact of our Theorem~\ref{thm:invariance} 
is that it gives the coin-tosser a chance to understand 
some of this infinite-dimensional theory. For example, 
note that for any fixed $u\ge 0$, 
the process $\{O^u_t\}_{t\ge 0}$ is an ordinary 
one-dimensional Ornstein--Uhlenbeck process. Therefore, 
Corollary~\ref{cor:genest-OU} can be stated, 
equivalently, as follows: 
\begin{corollary} 
   Let $E$ and $H$ be as in Theorem~\ref{thm:genest}. 
   Then there exists a finite 
   constant $A_{\ref{eq:co:genest}}>1$ such that for all 
   $z\ge 1$ and $u\ge 0$, 
   \begin{equation}\label{eq:co:genest} 
      A_{\ref{eq:co:genest}}^{-1} \K_E\left( 
      \frac{1}{z^2}\right) \bar\Phi(z) \le 
      \P\left\{ \sup_{t\in E} O^u_t \ge z 
      \right\} \le A_{\ref{eq:co:genest}} 
      \K_E\left( \frac{1}{z^2}\right) \bar\Phi(z). 
   \end{equation} 
\end{corollary} 
 
Similarly, the methods of this paper yield the following. 
We omit the details. 
\begin{corollary} 
   If $E$ and $H$ are as in Theorem~\ref{thm:inttest1}, 
   \begin{equation}\label{eq:multifractal}\begin{split} 
      &\sup_{t\in E} \limsup_{u\to\infty} 
         \left[ O^u_t - H(u)\sqrt{u} \right] > 0\ 
         \Longleftrightarrow\ \Psi_H(E)=+\infty\\ 
      &\Hdim\left\{ t\in[0,1]:\ 
         \limsup_{u\to\infty} \left[ 
         O^u_t  - H(u)\sqrt{u} \right] \ge 0\right\} 
         =\min\left( \frac{4-\d(H)}{2} ~,~ 1\right). 
   \end{split}\end{equation} 
\end{corollary} 
This is a multi-fractal extension of the main 
result of~\ocite{mountford} and extends some 
of the latter's infinite-dimensional potential theory. 
The results of this section seem to be new. 
 
\section{Concluding Remarks and Open Problems} 
 
The single-most important problem left open here is 
to remove the normality assumption in 
Theorems~\ref{thm:Hdim} and~\ref{thm:genest}. 
For instance, these theorems are not known to 
hold in the most important case where the increments 
are Rademacher variables. 
 
\begin{pbm} 
   Do Theorems~\ref{thm:Hdim} and~\ref{thm:genest} 
   hold for all incremental distributions $\nu$ that 
   have mean zero, variance one, and $2+\e$ moments 
   for some $\e>0$? 
\end{pbm} 
We suspect the answer is yes, but have no proof 
in any but the Gaussian case. 
 
As regards our invariance principles, we cannot 
resolve the following: 
 
\begin{pbm} 
   Does Theorem~\ref{thm:DS-REC} hold for $\e=0$? 
\end{pbm} 
We do not have a plausible conjecture in 
either direction. 
 
There is a large literature on tails of highly-oscillatory 
Gaussian random fields. See, for instance,~\ocite{pickands} and 
\ocite{qualls1}; see~\ocite{berman} for a pedagogic account 
as well as further references. In their simplest non-trivial 
setting, these works seek to find 
good asymptotic estimates for the tails of the distribution 
of $\sup_{t\in E}g(t)$ where $g$ is a stationary centered 
Gaussian random field that satisfies $\E\{ |g(0)-g(t)|^2 \}= 
1+c(1+o(1))|t|^\alpha$ as $|t|\to 0$. The ``time-set'' 
$E$ is often an interval or, more generally, a hyper-cube. 
What if $E$ is a fractal set? More generally, one can ask: 
\begin{pbm} 
   Do the results of \S\ref{sec:OU} have analogues for 
   more general Gaussian random fields? 
\end{pbm} 
 
There are a number of other interesting a.s.~properties 
of random walks one of which is the following due 
to~\ocite{chung}: Suppose $\{\xi_i\}_{i=1}^\infty$ are 
i.i.d.\@, mean-zero variance-one, and $\xi_1\in L^3(\P)$. 
Then $s_n=\xi_1+\cdots+\xi_n$ satisfies 
\begin{equation} 
   \liminf_{n\to\infty} \max_{1\le j\le n} 
   \frac{|s_j|}{\sqrt{n/\ln\ln n}} 
   =\frac{\pi}{\sqrt 8}\quad\text{ a.s.} 
\end{equation} 
\ocite{chung} contains the corresponding integral test. 
In the context of dynamical walks let us state, without 
proof, the following: If, in addition, $\xi_1\in L^4(\P)$, then 
\begin{equation} 
   \text{Chung's LIL is dynamically stable.} 
\end{equation} 
That is, with probability one, 
\begin{equation} 
   \liminf_{n\to\infty} \max_{1\le j\le n} 
   \frac{|S_j(t)|}{\sqrt{n/\ln\ln n}} 
   =\frac{\pi}{\sqrt 8}\qquad{}^\forall t\ge 0. 
\end{equation} 
\begin{pbm} 
   What can one say about the set of times $t\in[0,1]$ 
   at which $\{ S_n(t)\}_{n=1}^\infty$ is below 
   $\sqrt{n}/H(n)$ infinitely often? 
\end{pbm} 
This is related to finding sharp estimates for the 
``lower tail'' of $\max_{1\le j\le n}|S_j(t)|$. At 
this point, we have only partial results along these 
directions. For instance, when $\nu$ is standard normal, 
we can prove the existence of a constant $A$ such that 
for all compact $E\subseteq[0,1]$, 
\begin{equation}\label{chungchung} 
   \frac{e^{\pi^2/(8\e_n^2)}}{A\e_n^2} 
   \le \P\left\{ \inf_{t\in[0,1]} \max_{1\le j\le n} 
   |S_j(t)| \le \e_n\sqrt n\right\} \le 
   \frac{A e^{\pi^2/(8\e_n^2)}}{\e_n^6}, 
\end{equation} 
for any $(0,1)$-valued 
$\{\e_n\}_{n=1}^\infty$ that tends to zero 
and $\liminf_n n\e_n^8>\pi/\sqrt 2$. 
The solution to the preceding problem would 
require, invariably, a tightening of this bound. 
In a companion article~\cite{KLM2} 
we prove that the right-hand
side of~\eqref{chungchung} is tight for the continuum-limit
of dynamical walks. The said theorem
uses a second-order eigenvalue estimate
of~\ocite{SL} which is not yet available in the context
of dynamical random walks. Thus it is natural to end 
the paper with
the following open problem. 

\begin{pbm}
   Is the right-hand side of~\eqref{chungchung} 
   is sharp up to a multiplicative constant?
\end{pbm}

\end{document}